%% file: main.tex
\documentclass[a4paper]{article}

\usepackage[british]{babel}
\usepackage[utf8]{inputenc}
\usepackage{amsmath}
\usepackage{graphicx}
\usepackage[colorinlistoftodos]{todonotes}
\usepackage{hyperref}
\usepackage{multicol}
\usepackage{amsfonts}
\usepackage[british]{babel}
\usepackage{maplestd2e}
\usepackage[mathscr]{eucal}
\newtheorem{teo}{Theorem}
\newtheorem{defi}{Definition}

\title{Modeling Coral Reef Bleaching Recovery Through KCC-Theory}

    \author{Solange Rutz\footnote{Federal University of Pernambuco — UFPE, Mathematics Department, Recife, PE, Brazil} , \text{Rafael Cavalcanti}\footnote{Federal University of Pernambuco — UFPE, Mathematics Department, Recife, PE, Brazil} }
\date{}
\begin{document}
\maketitle
\begin{abstract}
We use Volterra-Hamilton systems theory and their associated cost functional to study the population dynamics and productive processes of coral reefs in recovery from bleaching and show that the cost of \mbox{production} remains the same after the process. The geometrical KCC-invariants are determined for the model proposed to describe the renewed symbiotic interaction between coral and algae.
\end{abstract}
\providecommand{\keywords}[1]
{
  \small	
  \textbf{\textit{Keywords---}} #1
}

\section{Introduction}
There is no doubt that increasing seawater temperature leads to coral bleaching {\cite{MS}}. This process occurs when corals are stressed by changes in environmental conditions such as temperature, light, or nutrients, leading to the expelling of the symbiotic algae which lives in the coral's tissues, causing it to  turn white. So global warming causes coral bleaching as increasing local seawater temperature stresses symbiotic algae (commonly called \textit{zooxanthellae}) in hermatypic coral (\textit{reef-building}) {\cite{PWG}} which leads to a breakdown in the symbiotic relationship between the coral animal and its zooxanthellae. This kind of seaweed has been severely affected by global warming around the world {\cite{ACB,CE}}. It is important to note that an individual coral is compound by thousand or even million of polyps which are animals of a few millimeters thick. The symbiotic algae living within the polyp makes energy from sunlight; they share that    energy with polyp in exchange for a comfortable environment and their interaction produce CaCO3 for the reef building. If there are no symbionts, the polyp run out of energy and dies within a few weeks or months, causing the coral to appear white or "bleached" {\cite{PW}}, unless they take more symbionts among those algae that naturally floating in the water around the coral barrier. Some coral reefs have been observed to recover from bleaching in appropriate situations {\cite{TR,PR}}. In order to model this recovering we suppose that before bleaching each polyp contains symbiotic algae living inside in a stable symbiotic relationship, and that there exist different kinds of algae outside the polyp, which benefit from the coral but do not influence it (\textit{commensal}) some of which are possibly better adapted to higher seawater temperatures.

\section{Interactive Population Dynamics}

Historically, a well-known population growth model is the Malthus model presented in 1798. This model considered that the rate of change of N is proportional to N according to a positive rate $\lambda>0$ called the specific growth rate, that is,
\begin{equation*}\label{malthus}
 \dfrac{dN}{dt} = \lambda N.
\end{equation*}
with solution
\begin{equation*}\label{solucaomalthus}
N(t) = N_0e^{\lambda t},    
\end{equation*}
for $N_0$ an initial condition. However, this model does not meet Hutchinson's third axiom.

Another well-known model is the logistic growth model which combines the exponential growth of N, with a reduction in growth that represents an environmental resistance. This model is given by
\begin{equation*}\label{logistica}
\dfrac{dN}{dt} = \lambda N - \dfrac{\lambda}{K} N^{2},    
\end{equation*}
with solution
\begin{equation*}\label{solucaologistica}
 N(t)=\dfrac{K}{1 + C_0 e^{-\lambda t}},   
\end{equation*}
where $C_0$ is a constant determined by the initial condition $N_0$ such that $C_0=\dfrac{K}{N_0} - 1$. In this case, the model meets all of Hutchinson's axioms.

The study of population growth dynamics can be extended to environments that have a greater number of species present, in this case, $\Sigma$ will represent a community, that is, a set of populations.

Assuming that such species do not interact with each other, we can use the logistic model to model the population dynamics of each species, so that
\begin{equation}\label{logisticanespecies}
\frac{dN^i}{dt} = \lambda_{(i)} N^i - \dfrac{\lambda_{(i)}}{K_{(i)}} (N^i)^2,
\end{equation}
where $N^i$ denotes the populational density of the species. The lack of interaction between species in (\ref{logisticanespecies}) is represented by the absence of $N^i N^j$ terms.

When there is interaction between species, the community is called simple community. For this case, (GAUSE, WITT) proposed
\begin{equation}
\left\{\begin{aligned}
\frac{dN^1}{dt} = \lambda_{(1)} N^1 - \dfrac{\lambda_{(1)}}{K_{(1)}} (N^1)^2 - \dfrac{\lambda_{(1)}\delta_{(1)}}{K_{(1)}} N^1 N^2 \\
\frac{dN^2}{dt} = \lambda_{(2)} N^2 - \dfrac{\lambda_{(2)}}{K_{(2)}} (N^2)^2 - \dfrac{\lambda_{(2)}\delta_{(2)}}{K_{(2)}} N^2 N^1,
\end{aligned}\right.
\end{equation}
where $\lambda_{(i)}, \delta_{(i)}, K_{(i)}$ are positive constants and the quantities $\delta_{(i)}$ represent how much the species $1$ is affected by the kind $2$ during the interaction.
\\
Let $\Sigma$ denote a population at a fixed location. Whenever we refer to $\Sigma$ it is implicit that there exists at least one individual living at this location. Let $N(t) \geq 0$ be the population density (or number of individuals per unit volume) in $\Sigma$ at time $t \in [0,T]$, $T>0$. We assume that $N(t)$ is a continuous function of time t
and satisfies Hutchinson's axioms {\cite{Hutchinson}} throughout this papper. We reserve the symbol $\Pi$ for the set of distinct species sharing the same location with no isolated species, i.e., it is a simple community. Let  $\Pi$ be a simple $n$ species community. We can describe this interaction by the system of equations:
\begin{equation}\label{eq1}
	\displaystyle\frac{dN^{i}}{dt}=\displaystyle\lambda_{(i)} N^{i}\left(1 - \frac{ N^{i}
	}{K_{(i)}} - \delta_{(j)}\frac{ N^{j}}{K_{(i)}}\right), \qquad \qquad i,j=1,2, \ldots, n.
\end{equation}
where $\lambda_{(i)}$ and $K_{(i)}$ are positive constants denoting \textit{intrinsic growth rate} and \textit{carrying capacity} for specie $i$, respectively. The coefficient $\delta_{(j)}$ represents how much the specie $i$ is affected by the specie $j$ in the interaction. The sign of $\delta_{(j)}$ tells what kind of interaction it is. Consider the system of two equations taking any $i,j$ $\in$ $\{1,\ldots,n\}$, $i \neq j$, in \eqref{eq1}. There are three possibilities to this system as follow:
\begin{itemize}
	\item Parasitism: $\delta_{(i)}> 0, \delta_{(j)} < 0$ or $\delta_{(i)}< 0, \delta_{(j)} > 0$;\\
	\item Competition: $\delta_{(i)}>0, \delta_{(j)} > 0$;\\
	\item Symbiosis: $\delta_{(i)}< 0, \delta_{(j)} < 0$.
\end{itemize}

\par According to {\cite{GW}} we see that Competition case is a Gause-Witt model with $n=2$. Then we have the following theorem (see {\cite{HV}}, page 21):
\begin{teo}\label{teo1}
 For a Gause-Witt model we have the following four cases: 
\begin{enumerate}
		
		\item If $\delta_{(1)}>\displaystyle\frac{K_{(1)}}{K_{(2)}}$ and $\delta_{(2)}>\displaystyle\frac{K_{(2)}}{K_{(1)}}$, then only one of the two species will persist after the competition and the winner will be determined entirely by the starting proportions.\\

		\item If $\delta_{(1)}>\displaystyle\frac{K_{(1)}}{K_{(2)}}$ and $\delta_{(2)}<\displaystyle\frac{K_{(2)}}{K_{(1)}}$, then the specie $1$ will be eliminated by the competition\\

		\item If $\delta_{((1)}<\displaystyle\frac{K_{(1)}}{K_{(2)}}$ and $\delta_{(2)}>\displaystyle\frac{K_{(2)}}{K_{(1)}}$, then the specie $2$ will be eliminated by the competition\\

		\item If $\delta_{(1)}<\displaystyle\frac{K_{(1)}}{K_{(2)}}$ and $\delta_{(2)}<\displaystyle\frac{K_{(2)}}{K_{(1)}}$, then both species persist together at equilibrium.
	\end{enumerate}	
\end{teo}

\par We are interested in study Gause-Witt equations whose coefficients satisfies hypothesis of items \textit{2 and 3}.\\
\par
\textit{Remark.} Gause-Witt model is just a particular case of a more general system of equation which describes ecological interaction of species in a simple community $\Pi$:
\begin{equation}\label{eq2}
	\displaystyle\frac{dN^{i}}{dt}=-\Gamma^{i}_{jk}N^{j}N^{k} + \lambda_{(i)}N^{i} + e^{i}, \qquad \quad i,j,k=1,\ldots,n
\end{equation}
where the $n^{3}$ quantities $\Gamma^{i}_{jk}$ are all constants (here, use is made of the Einstein summation convention on summing over repeated upper and lower indices with the only exception being $\lambda_{(i)}N^{i}$ where the parentheses indicates no summation).
\par Now we introduce a natural measure of production $x^{i}$ of a population $N^{i}(t)$, the Volterra's Production Variable {\cite{V}}, by defining
\begin{equation}\label{eq3}
	x^{i}(t)=k_{i}\displaystyle\int_{0}^{t}N^{i}(\tau)d\tau + x^{i}(0)
\end{equation}
where $k_{i}>0$ is the per capita production rate.
\par An n-species Volterra-Hamilton (VH) system $(\Pi, \Gamma)$ is pair consisting of $\Pi$, a set of n producer populations whose sizes are denoted by $N^{1}, \ldots, N^{n}$, together with a system of equations $\Gamma$ formed by \eqref{eq2} and \eqref{eq3}:
\begin{equation}\label{eq4}
	\Gamma:\begin{cases}
		\displaystyle\frac{dx^{i}}{dt}=k_{(i)}N^{i}  \\ 
		\displaystyle\frac{dN^{i}}{dt}=-\Gamma^{i}_{jk}N^{j}N^{k} + \lambda_{(i)}N^{i}+e^{i}
	\end{cases}
\end{equation}

\section{Bleaching Recovery Model}
Let $N^{1}(t)$, $N^{2}(t)$ and $N^{3}(t)$ be continuous functions of time which denote coral, symbiotic alga and commensal alga population density, respectively. We split this modelling in three stages: $(I)$ Commensal + Symbiosis; $(II)$ Symbiosis + Competition; $(III)$ Symbiosis. Here, it is initially assumed  $\lambda_{(N^{1})}=\lambda_{(N^{2})}=\lambda_{(N^{3})}=\lambda$, where these constants have the same meaning as in \eqref{eq1}. Following this symbiont assumptions, we can describe these three stages of interactions between coral and algae.\\

\par \textit{Remark.} The first stage describes how these three species live in the coral reef barrier before bleaching; Second stage is the dynamic produced by bleaching, but we will focus our attention at the competition between the algae; in the last one, we suppose alga $N^{3}$ will develop a symbiotic relation with the coral which was invaded by the outside algae, creating the condition to stop bleaching and start the recovering process.

\subsection{ Commensal + Symbiosis}
\par At this stage we suppose water temperature is adequate for both species of Algaes and to the Coral. First,  note that algae $N^{3}$ lives outside  the \textit{Polyp} (commensal relation), then this interaction is beneficial only one to alga $N^{3}$. By the other hand, alga $N^{2}$ and coral have a symbiotic interaction. So, we can describe this relation by extended Gause-Witt equations \eqref{eq1} as follow: 
\begin{center}
	\begin{equation}\label{eq5}
		\begin{cases}
			\displaystyle\frac{d N^{1}}{d t}=\lambda N^{1}-\frac{\lambda\left(N^{1}\right)^{2}}{K_{(1)}}+\delta_{(1)} \frac{\lambda N N^{2}}{K_{(1)}} \\ \\
			\displaystyle\frac{d N^{2}}{d t}=\lambda N^{2}-\frac{\lambda\left(N^{2}\right)^{2}}{K_{(2)}}+\delta_{(2)} \frac{\lambda N^{2} N^{1}}{K_{(2)}} \\ \\
			\displaystyle\frac{d N^{3}}{d t}=\lambda N^{3}-\frac{\lambda\left(N^{3}\right)^{2}}{K_{(3)}}+\delta_{(3)} \frac{\lambda N^{3} N^{1}}{K_{(3)}}
		\end{cases}
	\end{equation}
\end{center}
where $\delta's >0$ describes the symbiosis. If $N^{1}$ was affected by $N^{3}$, there would be a $4th$ term in the $1st$ equation .

\subsection{ Symbiosis + Competition}

\par \quad Here we assume that water warming is less lethal to $N^{3}$ than to $N^{2}$. The increasing water temperature produce a decreasing population density of algae $N^{2}$, since this specie is not adjusted to live in these conditions. This situation provides adequate conditions to algae $N^{3}$ penetrate the Polyp to establish the symbiotic relation that coral needs to live. In this case we have a dynamic where each specie interact to each other described as follow:
\begin{equation}\label{eq6}
	\begin{cases}
		\displaystyle\frac{d N^{1}}{d t}=\lambda N^{1}-\frac{\lambda\left(N^{1}\right)^{2}}{K_{(1)}}+\tilde{\delta}_{(1)} \frac{\lambda N^{1}\left(N^{2}+N^{3}\right)}{K_{(1)}}  \\ \\ 
		\displaystyle\frac{d N^{2}}{d t}=\lambda N^{2}-\frac{\lambda\left(N^{2}\right)^{2}}{K_{(2)}}+\delta_{(2)} \frac{\lambda N^{2} N^{1}}{K_{(2)}}-\mu_{(2)} \frac{\lambda N^{2} N^{3}}{K_{(2)}} \\ \\ \displaystyle\frac{d N^{3}}{d t}=\lambda N^{3}-\frac{\lambda\left(N^{3}\right)^{2}}{K_{(3)}}+\delta_{(3)} \frac{\lambda N^{3} N^{1}}{K_{(3)}}-\mu_{(3)} \frac{\lambda N^{2} N^{3}}{K_{(3)}}
	\end{cases}
\end{equation}

where $\mu_{(i)}$ ($i=2,3$) are positive contants and $\mu_{(i)}$ is the impact that $N^{i}$ suffers by interection with especie $N^{j}$, for $i, j \in\{2,3\}$. We expect that competition between $A^{1}$ and $A^{2}$ is so strong that we can assume $\mu_{(2)}, \mu_{(3)} \gg \tilde{\delta}_{(1)}, \delta_{(2)}, \delta_{(3)}$. Therefore \eqref{eq6} becomes a classical Gause-Witt competition system:
\begin{equation}\label{eq7}
	\begin{cases}
		\displaystyle\frac{d N^{2}}{d t}=\lambda N^{2}-\frac{\lambda\left(N^{2}\right)^{2}}{K_{(2)}}-\mu_{(2)} \frac{\lambda N^{2} N^{3}}{K_{(2)}} \\ \\ \displaystyle\frac{d N^{3}}{d t}=\lambda N^{3}-\frac{\lambda\left(N^{3}\right)^{2}}{K_{(3)}}-\mu_{(3)} \frac{\lambda N^{2} N^{3}}{K_{(3)}}
	\end{cases}
\end{equation}

\par As we have supposed that warmer water is more lethal to $N^{2}$ than $N^{3}$, then $\mu_{(3)}<\mu_{(2)}$ because this competition is harder $N^{2}$. Thus, by item 2 of theorem 1 we can conclude that $N^{2}$ is eliminated by the competition described in \eqref{eq7}.

\subsection{Symbiosis} 

\par After elimination of $N^{2}$ by competition with $N^{3}$, the coral $N^{1}$ has a new alga population to establish a symbiotic relation and then stop bleaching. The situation before bleaching and after recovering is quite the same in the sense of system of equations as follow: 
\begin{equation}\label{eq8}
	\begin{cases}
		\displaystyle\frac{d N^{1}}{d t}=\lambda N^{1}-\frac{\lambda\left(N^{1}\right)^{2}}{K_{(1)}}+\delta_{(1)} \frac{\lambda N^{1} N^{2}}{K_{(1)}}  \\ \\ 
		\displaystyle\frac{d N^{3}}{d t}=\lambda N^{3}-\frac{\lambda\left(N^{3}\right)^{2}}{K_{(3)}}+\delta_{(3)} \frac{\lambda N^{3} N^{1}}{K_{(3)}}
	\end{cases}
\end{equation}
\par \textit{Remark.} Equations \eqref{eq8} have the same form of the Gause-Witt system to describe interaction of $N^{1}$ and $N^{2}$ in \eqref{eq5}. This occurs because $N^{2}$ is supplanted by $N^{3}$.

\section{Proposal of the model}
Before bleaching disruption, it is known that coral and symbiotic alga develop a by-product as a result of their interaction. The same occurs after bleaching recovery since we are assuming alga $N^{3}$  becomes the symbiotic alga before the coral dies completely. Volterra-Hamilton is well suited to describe this production. For simplicity, we suppose all three populations have the same percapita rate of production(set $k_{(i)}=1, i=1,2,3$), so the production before bleaching is given by:
\begin{equation}\label{symbi11}
    \frac{d x^{1}}{d t}=N^{1}, \quad \frac{d N^{1}}{d t}=\lambda N^{1}-\frac{\lambda\left(N^{1}\right)^{2}}{K_{(1)}}+\delta_{(1)} \frac{\lambda N^{1} N^{2}}{K_{(1)}}
\end{equation}
and
\begin{equation}\label{symbi12}
    \frac{d x^{2}}{d t}=N^{2}, \quad \frac{d N^{2}}{d t}=\lambda N^{2}-\frac{\lambda\left(N^{2}\right)^{2}}{K_{(2)}}+\delta_{(2)} \frac{\lambda N^{2} N^{1}}{K_{(2)}}
\end{equation}
where the quantities $x^{1}$ and $x^{2}$ are the Volterra production variable \mbox{corresponding} to each specie.
\par Using the change of parameter $ds=\lambda e^{\lambda t}dt$ we obtain a system equivalent to equations \eqref{symbi11} and \eqref{symbi12}
\begin{equation}\label{eq19}
	\begin{cases}
		\displaystyle\frac{d^{2} x^{1}}{d s^{2}}+\frac{\lambda}{K_{(1)}}\left(\frac{d x^{1}}{d s}\right)^{2}-\frac{\lambda \delta_{(1)}}{K_{(1)}}\left(\frac{d x^{1}}{d s}\right)\left(\frac{d x^{2}}{d s}\right)=0 \\ \\
		\displaystyle\frac{d^{2} x^{2}}{d s^{2}}+\frac{\lambda}{K_{(2)}}\left(\frac{d x^{2}}{d s}\right)^{2}-\frac{\lambda \delta_{(2)}}{K_{(2)}}\left(\frac{d x^{2}}{d s}\right)\left(\frac{d x^{1}}{d s}\right)=0
	\end{cases}
\end{equation}
\par One can prove that
\begin{equation}\label{eq20}
	\displaystyle F(x, d x)=F\left(x^{1}, x^{2}, N^{1}, N^{2}\right)=e^{\psi\left(x^{1}, x^{2}\right)} \frac{\left(N^{2}\right)^{1+(1 / \lambda)}}{\left(N^{1}\right)^{1 / \lambda}}
\end{equation}
is conserved along the flow {(\ref{eq19})}, i.e., $dF/ds=0$, even these equations are not Euler-Lagrange for the functional $F$, where $\psi$ is of the form:
\begin{equation}\label{eq21}
	\psi\left(x^{1}, x^{2}\right)=A x^{1}+B x^{2},
\end{equation}
with
\[
A=-\displaystyle\left(\frac{\lambda \delta_{(2)}}{K_{(2)}}+\frac{K_{(2)}+\delta_{(2)} K_{(1)}}{K_{(1)} K_{(2)}}\right),\]
\[
B=\left(\frac{-\lambda \delta_{(1)}}{K_{(1)}}+\frac{(1+\lambda)\left(K_{(1)}+\delta_{(1)} K_{(2)}\right)}{K_{(1)} K_{(2)}}\right).
\]
\par By symmetry, the system that describes the dynamics after recovery is given by
\begin{equation}\label{eq22}
	\begin{cases}
		\displaystyle\frac{d^{2} x^{1}}{d s^{2}}+\frac{\lambda}{K_{(1)}}\left(\frac{d x^{1}}{d s}\right)^{2}-\frac{\lambda \delta_{(1)}}{K_{(1)}}\left(\frac{d x^{1}}{d s}\right)\left(\frac{d x^{3}}{d s}\right)=0 \\ \\
		\displaystyle\frac{d^{2} x^{3}}{d s^{2}}+\frac{\lambda}{K_{(3)}}\left(\frac{d x^{3}}{d s}\right)^{2}-\frac{\lambda \delta_{(3)}}{K_{(3)}}\left(\frac{d x^{1}}{d s}\right)\left(\frac{d x^{3}}{d s}\right)=0
	\end{cases}
\end{equation}
in an intrinsic time scale $s=e^{\lambda t}$, longer than $t$. Replacing $N^{2}$ for $N^{3}$, $x^{2}$ for $x^{3}$, $K_{(2)}$ for $K_{(3)}$ and $\delta_{(2)}$ for $\delta_{(3)}$ in \eqref{eq20} and \eqref{eq21} we conclude that after recovering, the cost of production is the same as before bleaching and $dF/ds=0$ along {(\ref{eq22})}, provided that we assume that the new alga replace the original one in the same ecological niche, or, in other words, $K_{(3)}$ apprimately equal $K_{(2)}$. This is interpreted as
representing an adaptation process, as oppose to an evolutionary one, where the cost of production is supposed to dimish, leading to a more efficient interaction pattern.
\par Now consider $\tilde{F}(x, d x)=F$, with $\psi(x)=A x^1+B x^2+\nu_3 x^1 x^2$. The Euler-Lagrange equation for $\tilde{F}$ were once obtained in \cite{SPor,Maple}, and they are:
\begin{equation}\label{metabolic}
	\begin{cases}
	\displaystyle\frac{d^2 x^1}{d s^2}-\lambda \nu_3 x^2\left(\frac{d x^1}{d s}\right)^2=0\\ \\
		\displaystyle\frac{d^2 x^2}{d s^2}+\frac{\lambda \nu_3}{\lambda+1} x^1\left(\frac{d x^2}{d s}\right)^2=0.
	\end{cases}
\end{equation}
The above system has some interesting properties and interpretations. First, note that the
functional F is conserved along trajectories of any of the three classical ecological
interactions, namely competition, symbiosis (or mutualism, which is a non-persistent
form of symbiosis) and parasitism, which differs just by the signs of their respective
interaction terms. The new, x-dependent system also preserves F along its trajectories,
but is the only one that is efficient, that is, satisfies Euler-Lagrange equations for F.
Being so, represents a metabolic, or x-dependent form of interactions, which may
indicate that any of the classical ecological systems actually conveys their interactions
through the exchange of products, which is particularly relevant as a model for
endosymbiosis, which results from an evolutionary process leading to a persistent
interdependence of the kind that occurs in the development of organs in complex, multi-
cellular species.

\section{KCC-Theory and Volterra-Hamilton System}
Let $(x^{1}, \ldots, x^{n})=(x)$, $\displaystyle\left(\frac{dx^{1}}{dt},\ldots, \frac{dx^{n}}{dt}\right)=\left(\frac{dx}{dt}\right)=(\dot{x})$ be $2n$ coordinates in an open connected subset $\Omega$ of the Euclidean $(2n)$-dimensional space $\mathbb{R}^{n} \times \mathbb{R}^{n}$. For our purpose, suppose that we have 
\begin{equation}\label{eq9}
	\frac{d^{2}x^{i}}{dt^{2}} + 2g^{i}(x,\dot{x})=0, \qquad i=1
	,\ldots, n,
\end{equation}
where each $g^{i}$ is $C^{\infty}$ in some neighborhoof of initial conditions $((x_{0}), (\dot{x}_{0}))\in \Omega$. The intrinsic geometry properties of {(\ref{eq8})} under non-singular tranformations of the type:
\begin{equation}\label{eq10}
	\begin{cases}
		\displaystyle\tilde{x}^{i}=f^{i}(x^{1}, \ldots, x^{n}), \qquad i=1, \ldots,n,\\ 
		\bar{t}=t
	\end{cases}
\end{equation}
are given by the five KCC-differential invariants, named after by D. Kosambi {\cite{Ko}}, E. Cartan \textcolor{blue}{\cite{Cart}} and S. S. Chern {\cite{Cher}}, given below. Let us first define the KCC-covariante differential of a contravariant vector field $\xi^{i}(x)$ on $\Omega$ by
\begin{equation}\label{eq11}
	\displaystyle\frac{\mathbb{D}\xi^{i}}{dt}=\frac{d\xi^{i}}{dt}+\frac{1}{2}g^{i}_{;r}\xi^{r}
\end{equation}
where the semi-colon indicates partial differntiation with respect to $\dot{x}^{r}$, and use of the Einstein summation convention on repeated indices. Using {($\ref{eq10}$)}, equation {(\ref{eq8})} becomes 
\begin{equation}\label{eq12}
	\displaystyle\frac{\mathbb{D}\xi^{i}}{dt}=\epsilon^{i}=\frac{1}{2}g^{i}_{;r}\dot{x}^{r}-g^{i}
\end{equation}
defining the first KCC-invariant of {(\ref{eq8})}, the contravariant vector field on $\Omega$, $\epsilon^{i}$, which represents an 'external force'. Varying trajectories $x^{i}(t)$ of {(\ref{eq8})} into nearby ones according to
\begin{equation}\label{eq13}
	\bar{x}^{i}(t)=x^{i}(t)+\xi^{i}(t)\eta
\end{equation}
where $\eta$ denotes a parameter, with $|\eta|$ small and $\xi^{i}(t)$ the components of some contravariant vector field defined along $x^{i}=x^{i}(t)$, we get, substituting {(\ref{eq12})} into {(\ref{eq8})} and taking the limit as $\eta \to 0$
\begin{equation}\label{eq14}
	\displaystyle\frac{d^{2}\xi^{i}}{dt^{2}}+g^{i}_{;r}\frac{d\xi^{r}}{dt}+g^{i}_{,r}\xi^{r}=0
\end{equation}\label{eq15}
where the comma indicates partial differentiation with respect to $x^{r}$. Using the KCC-covariant differentiation \eqref{eq10} we can express this as
\begin{equation}
	\displaystyle\frac{\mathbb{D}^{2}\xi^{i}}{dt^{2}}=\mathcal{P}^{i}_{r}\xi^{r},
\end{equation}
where 
\begin{equation}\label{eq16}
	\mathcal{P}^{i}_{j}=-g^{i}_{,j}-\frac{1}{2}g^{r}g^{i}_{;r;j}+\frac{1}{2}\dot{x}^{r}g^{i}_{,r;j}+\frac{1}{4}g^{i}_{;r}g^{r}_{;j}.
\end{equation}
The tensor $\mathcal{P}^{i}_{j}$ is the second KCC-invariant of \eqref{eq8}. The third, fourth and fifth invariants are:
\begin{equation}\label{eq17}
	\begin{cases}
		\mathcal{R}^{i}_{jk}=\displaystyle\frac{1}{3}(\mathcal{P}^{i}_{j;k}-\mathcal{P}^{i}_{k;j})\\
		\mathcal{B}^{i}_{jkl}=\mathcal{R}^{i}_{jk;l}\\
		\mathcal{D}^{i}_{jkl}=g^{i}_{;j;k;l}.
	\end{cases}
\end{equation}
The main result of KCC-theory is the following:
\begin{teo}[\cite{equiv}]\label{teo2}
	Two system of the form \eqref{eq8} on $\Omega$ are equivalent relative to \eqref{eq9} if and only if the five KCC-invariants are equivalent. In particular, there exist coordinates ($\bar{x}$) for which $g^{i}(\bar{x},\dot{\bar{x}},t)$ all vanish if and only if all KCC-invariants are zero. The tensor $\mathbb{D}$ vanishes if and only if $g^{i}$ is quadratic in ($\dot{x}$), in the case when the first KCC-invariant vanishes.
\end{teo}
\begin{defi}
Let  $\gamma(t)=\left(x^{i}(t)\right) \in U \subset \Omega$ be a path of \eqref{eq9}. If any other 
path with initial conditions close enough at $t=t_{0}$ remains close to $\gamma(t)$ for all $t>t_{0}$, we say that $\gamma(t)$ is a trajectory Jacobi stable. We define (3.2) to be Jacobi stable if all its solutions are Jacobi stable. Otherwise, we say that (3.2) is Jacobi unstable.
\end{defi}
\begin{teo}[\cite{Boehmer}]\label{teo3}
The trajectories of \eqref{eq9} are Jacobi stable if and only if the real part of the eigenvalues of the tensor $\mathcal{P}^{i}_{j}$ are strictly negative everywhere, and Jacobi unstable, otherwise. 
\end{teo}
\par Let us now introduce the notion of a n-dimensional Finsler space as a manifold where, given a coordinate system $(x)$ and a curve $x^{i}=x^{i}(t)$, the norm of a tangent vector $\dot{x}^{i}$ to the curve at each point $P$ on $x^{i}(t)$ is given by  the positive metric function $F$, $\left|\dot{x}^{i}\right|=F(x, \dot{x})$, where $F$ is positively homogeneous of degree 1 in $\dot{x}^{i}$. From $F$, a metric tensor is defined as $g_{i j}(x, \dot{x})=\displaystyle\frac{1}{2}\left(\partial^{2} F^{2} / \partial \dot{x}^{i} \partial \dot{x}^{j}\right)$, which must be regular in an open region of the tangent bundle, the collection of all tangent vectors to the manifold, and which excludes the origins. The use of the calculus of variations for $F$ leads to \eqref{eq9} with $g^{i}(x, \dot{x}, t)=\gamma_{j k}^{i}(x, \dot{x}) \dot{x}^{j} \dot{x}^{k}$, where the $\gamma_{j k}^{i}$ are the Levi-Cività symbols for the Finsler metric tensor $g_{i j}(x, \dot{x})$.
Berwald's Gaussian curvature $\mathcal{K}$ for two-dimensional Finsler spaces is defined from his famous formula \cite{ingarden} 
\begin{equation}\label{bgcK}
    \mathcal{R}_{j k}^{i}=F \mathcal{K} m^{i}\left(l_{j} m_{k}-l_{k} m_{j}\right)
\end{equation}
where $\mathscr{R}_{j k}^{i}$ is given by the first equation in \eqref{eq17}, $l^{i}=\dot{x}^{i} / F$ is the unit vector in the $\dot{x}^{i}$ direction, and $m^{i}$ the unique (up to orientation) unit vector perpendicular to $l^{i}$. Lowering the index on $m^{i}$ via the metric tensor gives $m_{i}$, which satisfies $F(x, m)=g_{i j}(x, \dot{x}) m^{i} m^{j} \equiv m_{i} m^{i}=1$. If our curvature $\mathcal{K}$ is bigger than zero everywhere, then trajectories oscillate back and forth, crossing the reference trajectory. In this case, we say \eqref{eq9} is Jacobi stable. If $\mathcal{K} \leq 0$ everywhere, trajectories diverge and system \eqref{eq9} is Jacobi unstable \cite{HV, diffusion}. This notion of stability is a Lyapunov notion, but it is a whole trajectory concept.
\begin{section}{Maple Computation}
Now, we will compute the five KCC-invariants of \eqref{eq19}, since it has the same form of \eqref{eq9}. For this computation we used the package \cite{SPor} at Maple \cite{Maple}.
\input{contas1}

One can realize that the system \eqref{eq22} is the same as \eqref{eq19} just by a replacement of $x^{2}$ by $x^{3}$. So, both systems are equivalent since their five KCC-invariants are equivalent. The production dynamic is not optimal before neither after the recovery of the coral reef as we can see using Theorem 3 on the signal eigenvalues of the second KCC-invariant of each case. Now, we will use the Berwald's Gaussian curvature \eqref{bgcK} to study the stability of the metabolic interaction given by \eqref{metabolic}:
\input{contas2}

\end{section}

\begin{section}{Conclusion and future work}

Our initial approach detailed in the section 3 is insufficient in terms of production between species. Although in the section 4 we consider that each specie has a natural product. With help of Maple computation we showed that the production described by equations \eqref{eq19} is not stable, using the Jacoby stability notion as in section 5. The same reasoning is used to \eqref{eq22}. The conclusion of this work is verify that symbiotic relation does not optimize production dynamics but the metabolic one does it.
\par A perspective of new researches from this work arise by the insertion of an external force in
the environment. If we consider again \eqref{symbi11} and \eqref{symbi12}, but now with the force $e^i=-\left(\delta_j^i \sigma_k(x)\right) N^j N^k$, and $\sigma_k(x)$ being g a smooth covariant vector field on production space there is a new dynamic
to investigate. 

\end{section}
\section*{Conflict of Interest}
The authors declare no conflicts of interest

\section*{Acknowledgements}
The second author was supported by CNPq.

\end{document}

%% file: contas1.tex
\pagestyle{empty}
\DefineParaStyle{Maple Heading 1}
\DefineParaStyle{Maple Text Output}
\DefineParaStyle{Maple Dash Item}
\DefineParaStyle{Maple Bullet Item}
\DefineParaStyle{Maple Normal}
\DefineParaStyle{Maple Heading 4}
\DefineParaStyle{Maple Heading 3}
\DefineParaStyle{Maple Heading 2}
\DefineParaStyle{Maple Warning}
\DefineParaStyle{Maple Title}
\DefineParaStyle{Maple Error}
\DefineCharStyle{Maple Hyperlink}
\DefineCharStyle{Maple 2D Math}
\DefineCharStyle{Maple Maple Input}
\DefineCharStyle{Maple 2D Output}
\DefineCharStyle{Maple 2D Input}
\begin{maplegroup}
\mapleresult
\begin{maplelatex}
\end{maplelatex}
\end{maplegroup}
\textbf{}\mapleinline{inert}{2d}{with(Finsler); 1}{\[\displaystyle \]}
\begin{maplegroup}
\mapleresult
\begin{maplelatex}
\mapleinline{inert}{2d}{[Dcoordinates, Hdiff, K, connection, init, metricfunction, tddiff]}{\[\displaystyle [{\it Dcoordinates},{\it Hdiff},K,{\it connection}\\
\mbox{},{\it init},{\it metricfunction},{\it tddiff}]\]}
\end{maplelatex}
\end{maplegroup}
\mapleinline{inert}{2d}{dimension := 2; 1}{\[\displaystyle \]}
\begin{maplegroup}
\mapleresult
\begin{maplelatex}
\mapleinline{inert}{2d}{2}{\[\displaystyle 2\]}
\end{maplelatex}
\end{maplegroup}
\begin{maplegroup}
\begin{mapleinput}
\mapleinline{active}{2d}{coordinates(x1, x2); 1}{\[\]}
\end{mapleinput}
\mapleresult
\begin{maplelatex}
\mapleinline{inert}{2d}{The coordinates are:}{\[\displaystyle \mbox {{\tt The coordinates are:}}\]}
\end{maplelatex}
\mapleresult
\begin{maplelatex}
\mapleinline{inert}{2d}{X^1 = x1}{\[\displaystyle X{\mbox {{\tt }}}^{\mbox {{\tt 1}}}={\it x1}\]}
\end{maplelatex}
\mapleresult
\begin{maplelatex}
\mapleinline{inert}{2d}{X^2 = x2}{\[\displaystyle X{\mbox {{\tt }}}^{2}={\it x2}\]}
\end{maplelatex}
\end{maplegroup}
\begin{maplegroup}
\begin{mapleinput}
\mapleinline{active}{2d}{Dcoordinates(N1, N2); 1}{\[\]}
\end{mapleinput}
\mapleresult
\begin{maplelatex}
\mapleinline{inert}{2d}{The d-coordinates are:}{\[\displaystyle \mbox {{\tt The d-coordinates are:}}\]}
\end{maplelatex}
\mapleresult
\begin{maplelatex}
\mapleinline{inert}{2d}{Y^1 = N1}{\[\displaystyle Y{\mbox {{\tt}}}^{\mbox {{\tt 1}}}={\it N1}\]}
\end{maplelatex}
\mapleresult
\begin{maplelatex}
\mapleinline{inert}{2d}{Y^2 = N2}{\[\displaystyle Y{\mbox {{\tt }}}^{2}={\it N2}\]}
\end{maplelatex}
\end{maplegroup}
\begin{maplegroup}
\begin{mapleinput}
\mapleinline{active}{2d}{G1 := (1/2)*lambda*N1/K1*N1-(1/2)*lambda*delta1*N1*N2/K1; 1}{\[\]}
\end{mapleinput}
\mapleresult
\begin{maplelatex}
\mapleinline{inert}{2d}{G1 := (1/2)*lambda*N1^2/K1-(1/2)*lambda*delta1*N1*N2/K1}{\[\displaystyle {\it G1}\, := \,1/2\,{\frac {\lambda\,{{\it N1}}^{2}}{{\it K1}}}-1/2\,{\frac {\lambda\,{\it delta1}\,{\it N1}\,{\it N2}}{{\it K1}}}\]}
\end{maplelatex}
\end{maplegroup}
\begin{maplegroup}
\begin{mapleinput}
\mapleinline{active}{2d}{G2 := (1/2)*lambda*N2/K2*N2-(1/2)*lambda*delta2*N1*N2/K2; 1}{\[\]}
\end{mapleinput}
\mapleresult
\begin{maplelatex}
\mapleinline{inert}{2d}{G2 := (1/2)*lambda*N2^2/K2-(1/2)*lambda*delta2*N1*N2/K2}{\[\displaystyle {\it G2}\, := \,1/2\,{\frac {\lambda\,{{\it N2}}^{2}}{{\it K2}}}-1/2\,{\frac {\lambda\,{\it delta2}\,{\it N1}\,{\it N2}}{{\it K2}}}\]}
\end{maplelatex}
\end{maplegroup}
\begin{maplegroup}
\begin{mapleinput}
\mapleinline{active}{2d}{connection(G1, G2); 1}{\[\]}
\end{mapleinput}
\mapleresult
\begin{maplelatex}
\mapleinline{inert}{2d}{G^x1 = (1/2)*lambda*N1^2/K1-(1/2)*lambda*delta1*N1*N2/K1}{\[\displaystyle G{\mbox {{\tt}}}^{{\it x1}}:=1/2\,{\frac {\lambda\,{{\it N1}}^{2}}{{\it K1}}}\\
\mbox{}-1/2\,{\frac {\lambda\,{\it delta1}\,{\it N1}\,{\it N2}}{{\it K1}}}\]}
\end{maplelatex}
\mapleresult
\begin{maplelatex}
\mapleinline{inert}{2d}{G^x2 = (1/2)*lambda*N2^2/K2-(1/2)*lambda*delta2*N1*N2/K2}{\[\displaystyle G{\mbox {{\tt}}}^{{\it x2}}:=1/2\,{\frac {\lambda\,{{\it N2}}^{2}}{{\it K2}}}\\
\mbox{}-1/2\,{\frac {\lambda\,{\it delta2}\,{\it N1}\,{\it N2}}{{\it K2}}}\]}
\end{maplelatex}
\end{maplegroup}
\begin{maplegroup}
\begin{mapleinput}
\mapleinline{active}{2d}{definetensor(epsilon[i] = N[i, -j]*Y[j]-2*G[i]); 1}{\[\]}
\end{mapleinput}
\begin{maplelatex}
\mapleinline{inert}{2d}{First KCC-invariant:}{\[\displaystyle \mbox {{\tt First KCC-invariant:}}\]}
\end{maplelatex}
\mapleresult
\begin{maplelatex}
\mapleinline{inert}{2d}{Typesetting:-mrow(Typesetting:-mrow(Typesetting:-mi("&varepsilon;", italic = "false", mathvariant = "normal"), Typesetting:-mo("&InvisibleTimes;", mathvariant = "normal", fence = "false", separator = "false", stretchy = "false", symmetric = "false", largeop = "false", movablelimits = "false", accent = "false", lspace = "0.0em", rspace = "0.0em"), Typesetting:-msup(Typesetting:-mi(" ", italic = "true", mathvariant = "italic"), Typesetting:-mi("i", italic = "true", mathvariant = "italic"), superscriptshift = "0")), Typesetting:-mo("=", mathvariant = "normal", fence = "false", separator = "false", stretchy = "false", symmetric = "false", largeop = "false", movablelimits = "false", accent = "false", lspace = "0.2777778em", rspace = "0.2777778em"), Typesetting:-mrow(Typesetting:-mrow(Typesetting:-mi("N", italic = "true", mathvariant = "italic"), Typesetting:-mo("&InvisibleTimes;", mathvariant = "normal", fence = "false", separator = "false", stretchy = "false", symmetric = "false", largeop = "false", movablelimits = "false", accent = "false", lspace = "0.0em", rspace = "0.0em"), Typesetting:-msup(Typesetting:-mi(" ", italic = "true", mathvariant = "italic"), Typesetting:-mi("i", italic = "true", mathvariant = "italic"), superscriptshift = "0"), Typesetting:-mo("&InvisibleTimes;", mathvariant = "normal", fence = "false", separator = "false", stretchy = "false", symmetric = "false", largeop = "false", movablelimits = "false", accent = "false", lspace = "0.0em", rspace = "0.0em"), Typesetting:-msub(Typesetting:-mi("", italic = "true", mathvariant = "italic"), Typesetting:-mrow(Typesetting:-mi("j", italic = "true", mathvariant = "italic")), subscriptshift = "0"), Typesetting:-mo("&InvisibleTimes;", mathvariant = "normal", fence = "false", separator = "false", stretchy = "false", symmetric = "false", largeop = "false", movablelimits = "false", accent = "false", lspace = "0.0em", rspace = "0.0em"), Typesetting:-mi("Y", italic = "true", mathvariant = "italic"), Typesetting:-mo("&InvisibleTimes;", mathvariant = "normal", fence = "false", separator = "false", stretchy = "false", symmetric = "false", largeop = "false", movablelimits = "false", accent = "false", lspace = "0.0em", rspace = "0.0em"), Typesetting:-msup(Typesetting:-mi(" ", italic = "true", mathvariant = "italic"), Typesetting:-mi("j", italic = "true", mathvariant = "italic"), superscriptshift = "0")), Typesetting:-mo("&minus;", mathvariant = "normal", fence = "false", separator = "false", stretchy = "false", symmetric = "false", largeop = "false", movablelimits = "false", accent = "false", lspace = "0.2222222em", rspace = "0.2222222em"), Typesetting:-mrow(Typesetting:-mn("2", mathvariant = "normal"), Typesetting:-mo("&InvisibleTimes;", mathvariant = "normal", fence = "false", separator = "false", stretchy = "false", symmetric = "false", largeop = "false", movablelimits = "false", accent = "false", lspace = "0.0em", rspace = "0.0em"), Typesetting:-mi("G", italic = "true", mathvariant = "italic"), Typesetting:-mo("&InvisibleTimes;", mathvariant = "normal", fence = "false", separator = "false", stretchy = "false", symmetric = "false", largeop = "false", movablelimits = "false", accent = "false", lspace = "0.0em", rspace = "0.0em"), Typesetting:-msup(Typesetting:-mi(" ", italic = "true", mathvariant = "italic"), Typesetting:-mi("i", italic = "true", mathvariant = "italic"), superscriptshift = "0"))))}{\[\displaystyle \varepsilon ~\mbox {{\tt }} ^{i }=N ~\mbox {{\tt }} ^{i }~_{j }~Y ~\mbox {{\tt }} ^{j }-2~G ~\mbox {{\tt }} ^{i }\]}
\end{maplelatex}
\end{maplegroup}
\begin{maplegroup}
\begin{mapleinput}
\mapleinline{active}{2d}{evalt(epsilon[i]); 1}{\[\]}
\end{mapleinput}
\mapleresult
\begin{maplelatex}
\mapleinline{inert}{2d}{epsilon*^i = 0}{\[\displaystyle \epsilon\,{\mbox {{\tt }}}^{i}=0\]}
\end{maplelatex}
\end{maplegroup}
\begin{maplegroup}
\begin{mapleinput}
\mapleinline{active}{2d}{show(B[i, -j]); 1}{\[\]}
\end{mapleinput}
\begin{maplelatex}
\mapleinline{inert}{2d}{Second KCC-invariant:}{\[\displaystyle \mbox {{\tt Second KCC-invariant:}}\]}
\end{maplelatex}
\mapleresult
\begin{maplelatex}
\mapleinline{inert}{2d}{Typesetting:-mrow(Typesetting:-mrow(Typesetting:-mi("B", italic = "true", mathvariant = "italic"), Typesetting:-mo("&InvisibleTimes;", mathvariant = "normal", fence = "false", separator = "false", stretchy = "false", symmetric = "false", largeop = "false", movablelimits = "false", accent = "false", lspace = "0.0em", rspace = "0.0em"), Typesetting:-msup(Typesetting:-mi(" ", italic = "true", mathvariant = "italic"), Typesetting:-mi("x1", italic = "true", mathvariant = "italic"), superscriptshift = "0"), Typesetting:-mo("&InvisibleTimes;", mathvariant = "normal", fence = "false", separator = "false", stretchy = "false", symmetric = "false", largeop = "false", movablelimits = "false", accent = "false", lspace = "0.0em", rspace = "0.0em"), Typesetting:-msub(Typesetting:-mi("", italic = "true", mathvariant = "italic"), Typesetting:-mrow(Typesetting:-mi("x1", italic = "true", mathvariant = "italic")), subscriptshift = "0")), Typesetting:-mo("=", mathvariant = "normal", fence = "false", separator = "false", stretchy = "false", symmetric = "false", largeop = "false", movablelimits = "false", accent = "false", lspace = "0.2777778em", rspace = "0.2777778em"), Typesetting:-mrow(Typesetting:-mfrac(Typesetting:-mn("1", mathvariant = "normal"), Typesetting:-mn("4", mathvariant = "normal"), linethickness = "1", denomalign = "center", numalign = "center", bevelled = "false"), Typesetting:-mo("&InvisibleTimes;", mathvariant = "normal", fence = "false", separator = "false", stretchy = "false", symmetric = "false", largeop = "false", movablelimits = "false", accent = "false", lspace = "0.0em", rspace = "0.0em"), Typesetting:-mfrac(Typesetting:-mrow(Typesetting:-msup(Typesetting:-mi("&lambda;", italic = "false", mathvariant = "normal"), Typesetting:-mn("2", mathvariant = "normal"), superscriptshift = "0"), Typesetting:-mo("&InvisibleTimes;", mathvariant = "normal", fence = "false", separator = "false", stretchy = "false", symmetric = "false", largeop = "false", movablelimits = "false", accent = "false", lspace = "0.0em", rspace = "0.0em"), Typesetting:-mi("&delta;1", italic = "true", mathvariant = "italic"), Typesetting:-mo("&InvisibleTimes;", mathvariant = "normal", fence = "false", separator = "false", stretchy = "false", symmetric = "false", largeop = "false", movablelimits = "false", accent = "false", lspace = "0.0em", rspace = "0.0em"), Typesetting:-mi("N2", italic = "true", mathvariant = "italic"), Typesetting:-mo("&InvisibleTimes;", mathvariant = "normal", fence = "false", separator = "false", stretchy = "false", symmetric = "false", largeop = "false", movablelimits = "false", accent = "false", lspace = "0.0em", rspace = "0.0em"), Typesetting:-mfenced(Typesetting:-mrow(Typesetting:-mrow(Typesetting:-mi("K1", italic = "true", mathvariant = "italic"), Typesetting:-mo("&InvisibleTimes;", mathvariant = "normal", fence = "false", separator = "false", stretchy = "false", symmetric = "false", largeop = "false", movablelimits = "false", accent = "false", lspace = "0.0em", rspace = "0.0em"), Typesetting:-mi("N1", italic = "true", mathvariant = "italic"), Typesetting:-mo("&InvisibleTimes;", mathvariant = "normal", fence = "false", separator = "false", stretchy = "false", symmetric = "false", largeop = "false", movablelimits = "false", accent = "false", lspace = "0.0em", rspace = "0.0em"), Typesetting:-mi("&delta;2", italic = "true", mathvariant = "italic")), Typesetting:-mo("&minus;", mathvariant = "normal", fence = "false", separator = "false", stretchy = "false", symmetric = "false", largeop = "false", movablelimits = "false", accent = "false", lspace = "0.2222222em", rspace = "0.2222222em"), Typesetting:-mrow(Typesetting:-mi("K2", italic = "true", mathvariant = "italic"), Typesetting:-mo("&InvisibleTimes;", mathvariant = "normal", fence = "false", separator = "false", stretchy = "false", symmetric = "false", largeop = "false", movablelimits = "false", accent = "false", lspace = "0.0em", rspace = "0.0em"), Typesetting:-mi("N2", italic = "true", mathvariant = "italic"), Typesetting:-mo("&InvisibleTimes;", mathvariant = "normal", fence = "false", separator = "false", stretchy = "false", symmetric = "false", largeop = "false", movablelimits = "false", accent = "false", lspace = "0.0em", rspace = "0.0em"), Typesetting:-mi("&delta;1", italic = "true", mathvariant = "italic")), Typesetting:-mo("&minus;", mathvariant = "normal", fence = "false", separator = "false", stretchy = "false", symmetric = "false", largeop = "false", movablelimits = "false", accent = "false", lspace = "0.2222222em", rspace = "0.2222222em"), Typesetting:-mrow(Typesetting:-mn("2", mathvariant = "normal"), Typesetting:-mo("&InvisibleTimes;", mathvariant = "normal", fence = "false", separator = "false", stretchy = "false", symmetric = "false", largeop = "false", movablelimits = "false", accent = "false", lspace = "0.0em", rspace = "0.0em"), Typesetting:-mi("K1", italic = "true", mathvariant = "italic"), Typesetting:-mo("&InvisibleTimes;", mathvariant = "normal", fence = "false", separator = "false", stretchy = "false", symmetric = "false", largeop = "false", movablelimits = "false", accent = "false", lspace = "0.0em", rspace = "0.0em"), Typesetting:-mi("N2", italic = "true", mathvariant = "italic"))), mathvariant = "normal")), Typesetting:-mrow(Typesetting:-msup(Typesetting:-mi("K1", italic = "true", mathvariant = "italic"), Typesetting:-mn("2", mathvariant = "normal"), superscriptshift = "0"), Typesetting:-mo("&InvisibleTimes;", mathvariant = "normal", fence = "false", separator = "false", stretchy = "false", symmetric = "false", largeop = "false", movablelimits = "false", accent = "false", lspace = "0.0em", rspace = "0.0em"), Typesetting:-mi("K2", italic = "true", mathvariant = "italic")), linethickness = "1", denomalign = "center", numalign = "center", bevelled = "false")))}
{\[\displaystyle B ~\mbox {{\tt }} ^{{\it x1} }~_{{\it x1} }=\frac{1}{4}~\frac{\lambda ^{2}~\delta \mbox {{\tt 1}}\\
\mbox{} ~{\it N2} ~\left({\it K1} ~{\it N1} ~\delta \mbox {{\tt 2}}\\
\mbox{} -{\it K2} ~{\it N2} ~\delta \mbox {{\tt 1}}\\
\mbox{} -2~{\it K1} ~{\it N2} \right)}{{\it K1} ^{2}~{\it K2} }\]}
\end{maplelatex}
\mapleresult
\begin{maplelatex}
\mapleinline{inert}{2d}{Typesetting:-mrow(Typesetting:-mrow(Typesetting:-mi("B", italic = "true", mathvariant = "italic"), Typesetting:-mo("&InvisibleTimes;", mathvariant = "normal", fence = "false", separator = "false", stretchy = "false", symmetric = "false", largeop = "false", movablelimits = "false", accent = "false", lspace = "0.0em", rspace = "0.0em"), Typesetting:-msup(Typesetting:-mi(" ", italic = "true", mathvariant = "italic"), Typesetting:-mi("x1", italic = "true", mathvariant = "italic"), superscriptshift = "0"), Typesetting:-mo("&InvisibleTimes;", mathvariant = "normal", fence = "false", separator = "false", stretchy = "false", symmetric = "false", largeop = "false", movablelimits = "false", accent = "false", lspace = "0.0em", rspace = "0.0em"), Typesetting:-msub(Typesetting:-mi("", italic = "true", mathvariant = "italic"), Typesetting:-mrow(Typesetting:-mi("x2", italic = "true", mathvariant = "italic")), subscriptshift = "0")), Typesetting:-mo("=", mathvariant = "normal", fence = "false", separator = "false", stretchy = "false", symmetric = "false", largeop = "false", movablelimits = "false", accent = "false", lspace = "0.2777778em", rspace = "0.2777778em"), Typesetting:-mrow(Typesetting:-mo("&uminus0;", mathvariant = "normal", fence = "false", separator = "false", stretchy = "false", symmetric = "false", largeop = "false", movablelimits = "false", accent = "false", lspace = "0.2222222em", rspace = "0.2222222em"), Typesetting:-mrow(Typesetting:-mfrac(Typesetting:-mn("1", mathvariant = "normal"), Typesetting:-mn("4", mathvariant = "normal"), linethickness = "1", denomalign = "center", numalign = "center", bevelled = "false"), Typesetting:-mo("&InvisibleTimes;", mathvariant = "normal", fence = "false", separator = "false", stretchy = "false", symmetric = "false", largeop = "false", movablelimits = "false", accent = "false", lspace = "0.0em", rspace = "0.0em"), Typesetting:-mfrac(Typesetting:-mrow(Typesetting:-msup(Typesetting:-mi("&lambda;", italic = "false", mathvariant = "normal"), Typesetting:-mn("2", mathvariant = "normal"), superscriptshift = "0"), Typesetting:-mo("&InvisibleTimes;", mathvariant = "normal", fence = "false", separator = "false", stretchy = "false", symmetric = "false", largeop = "false", movablelimits = "false", accent = "false", lspace = "0.0em", rspace = "0.0em"), Typesetting:-mi("&delta;1", italic = "true", mathvariant = "italic"), Typesetting:-mo("&InvisibleTimes;", mathvariant = "normal", fence = "false", separator = "false", stretchy = "false", symmetric = "false", largeop = "false", movablelimits = "false", accent = "false", lspace = "0.0em", rspace = "0.0em"), Typesetting:-mi("N1", italic = "true", mathvariant = "italic"), Typesetting:-mo("&InvisibleTimes;", mathvariant = "normal", fence = "false", separator = "false", stretchy = "false", symmetric = "false", largeop = "false", movablelimits = "false", accent = "false", lspace = "0.0em", rspace = "0.0em"), Typesetting:-mfenced(Typesetting:-mrow(Typesetting:-mrow(Typesetting:-mi("K1", italic = "true", mathvariant = "italic"), Typesetting:-mo("&InvisibleTimes;", mathvariant = "normal", fence = "false", separator = "false", stretchy = "false", symmetric = "false", largeop = "false", movablelimits = "false", accent = "false", lspace = "0.0em", rspace = "0.0em"), Typesetting:-mi("N1", italic = "true", mathvariant = "italic"), Typesetting:-mo("&InvisibleTimes;", mathvariant = "normal", fence = "false", separator = "false", stretchy = "false", symmetric = "false", largeop = "false", movablelimits = "false", accent = "false", lspace = "0.0em", rspace = "0.0em"), Typesetting:-mi("&delta;2", italic = "true", mathvariant = "italic")), Typesetting:-mo("&minus;", mathvariant = "normal", fence = "false", separator = "false", stretchy = "false", symmetric = "false", largeop = "false", movablelimits = "false", accent = "false", lspace = "0.2222222em", rspace = "0.2222222em"), Typesetting:-mrow(Typesetting:-mi("K2", italic = "true", mathvariant = "italic"), Typesetting:-mo("&InvisibleTimes;", mathvariant = "normal", fence = "false", separator = "false", stretchy = "false", symmetric = "false", largeop = "false", movablelimits = "false", accent = "false", lspace = "0.0em", rspace = "0.0em"), Typesetting:-mi("N2", italic = "true", mathvariant = "italic"), Typesetting:-mo("&InvisibleTimes;", mathvariant = "normal", fence = "false", separator = "false", stretchy = "false", symmetric = "false", largeop = "false", movablelimits = "false", accent = "false", lspace = "0.0em", rspace = "0.0em"), Typesetting:-mi("&delta;1", italic = "true", mathvariant = "italic")), Typesetting:-mo("&minus;", mathvariant = "normal", fence = "false", separator = "false", stretchy = "false", symmetric = "false", largeop = "false", movablelimits = "false", accent = "false", lspace = "0.2222222em", rspace = "0.2222222em"), Typesetting:-mrow(Typesetting:-mn("2", mathvariant = "normal"), Typesetting:-mo("&InvisibleTimes;", mathvariant = "normal", fence = "false", separator = "false", stretchy = "false", symmetric = "false", largeop = "false", movablelimits = "false", accent = "false", lspace = "0.0em", rspace = "0.0em"), Typesetting:-mi("K1", italic = "true", mathvariant = "italic"), Typesetting:-mo("&InvisibleTimes;", mathvariant = "normal", fence = "false", separator = "false", stretchy = "false", symmetric = "false", largeop = "false", movablelimits = "false", accent = "false", lspace = "0.0em", rspace = "0.0em"), Typesetting:-mi("N2", italic = "true", mathvariant = "italic"))), mathvariant = "normal")), Typesetting:-mrow(Typesetting:-msup(Typesetting:-mi("K1", italic = "true", mathvariant = "italic"), Typesetting:-mn("2", mathvariant = "normal"), superscriptshift = "0"), Typesetting:-mo("&InvisibleTimes;", mathvariant = "normal", fence = "false", separator = "false", stretchy = "false", symmetric = "false", largeop = "false", movablelimits = "false", accent = "false", lspace = "0.0em", rspace = "0.0em"), Typesetting:-mi("K2", italic = "true", mathvariant = "italic")), linethickness = "1", denomalign = "center", numalign = "center", bevelled = "false"))))}{\[\displaystyle B ~\mbox {{\tt }} ^{{\it x1} }~_{{\it x2} }=-\frac{1}{4}~\frac{\lambda ^{2}~\delta \mbox {{\tt 1}}\\
\mbox{} ~{\it N1} ~\left({\it K1}\\
\mbox{} ~{\it N1} ~\delta \mbox {{\tt 2}} -{\it K2}\\
\mbox{} ~{\it N2} ~\delta \mbox {{\tt 1}}\\
\mbox{} -2~{\it K1}\\
\mbox{} ~{\it N2} \right)}{{\it K1}\\
\mbox{} ^{2}~{\it K2}\\
\mbox{} }\]}
\end{maplelatex}
\mapleresult
\begin{maplelatex}
\mapleinline{inert}{2d}{Typesetting:-mrow(Typesetting:-mrow(Typesetting:-mi("B", italic = "true", mathvariant = "italic"), Typesetting:-mo("&InvisibleTimes;", mathvariant = "normal", fence = "false", separator = "false", stretchy = "false", symmetric = "false", largeop = "false", movablelimits = "false", accent = "false", lspace = "0.0em", rspace = "0.0em"), Typesetting:-msup(Typesetting:-mi(" ", italic = "true", mathvariant = "italic"), Typesetting:-mi("x2", italic = "true", mathvariant = "italic"), superscriptshift = "0"), Typesetting:-mo("&InvisibleTimes;", mathvariant = "normal", fence = "false", separator = "false", stretchy = "false", symmetric = "false", largeop = "false", movablelimits = "false", accent = "false", lspace = "0.0em", rspace = "0.0em"), Typesetting:-msub(Typesetting:-mi("", italic = "true", mathvariant = "italic"), Typesetting:-mrow(Typesetting:-mi("x1", italic = "true", mathvariant = "italic")), subscriptshift = "0")), Typesetting:-mo("=", mathvariant = "normal", fence = "false", separator = "false", stretchy = "false", symmetric = "false", largeop = "false", movablelimits = "false", accent = "false", lspace = "0.2777778em", rspace = "0.2777778em"), Typesetting:-mrow(Typesetting:-mfrac(Typesetting:-mn("1", mathvariant = "normal"), Typesetting:-mn("4", mathvariant = "normal"), linethickness = "1", denomalign = "center", numalign = "center", bevelled = "false"), Typesetting:-mo("&InvisibleTimes;", mathvariant = "normal", fence = "false", separator = "false", stretchy = "false", symmetric = "false", largeop = "false", movablelimits = "false", accent = "false", lspace = "0.0em", rspace = "0.0em"), Typesetting:-mfrac(Typesetting:-mrow(Typesetting:-msup(Typesetting:-mi("&lambda;", italic = "false", mathvariant = "normal"), Typesetting:-mn("2", mathvariant = "normal"), superscriptshift = "0"), Typesetting:-mo("&InvisibleTimes;", mathvariant = "normal", fence = "false", separator = "false", stretchy = "false", symmetric = "false", largeop = "false", movablelimits = "false", accent = "false", lspace = "0.0em", rspace = "0.0em"), Typesetting:-mi("&delta;2", italic = "true", mathvariant = "italic"), Typesetting:-mo("&InvisibleTimes;", mathvariant = "normal", fence = "false", separator = "false", stretchy = "false", symmetric = "false", largeop = "false", movablelimits = "false", accent = "false", lspace = "0.0em", rspace = "0.0em"), Typesetting:-mi("N2", italic = "true", mathvariant = "italic"), Typesetting:-mo("&InvisibleTimes;", mathvariant = "normal", fence = "false", separator = "false", stretchy = "false", symmetric = "false", largeop = "false", movablelimits = "false", accent = "false", lspace = "0.0em", rspace = "0.0em"), Typesetting:-mfenced(Typesetting:-mrow(Typesetting:-mrow(Typesetting:-mi("K1", italic = "true", mathvariant = "italic"), Typesetting:-mo("&InvisibleTimes;", mathvariant = "normal", fence = "false", separator = "false", stretchy = "false", symmetric = "false", largeop = "false", movablelimits = "false", accent = "false", lspace = "0.0em", rspace = "0.0em"), Typesetting:-mi("N1", italic = "true", mathvariant = "italic"), Typesetting:-mo("&InvisibleTimes;", mathvariant = "normal", fence = "false", separator = "false", stretchy = "false", symmetric = "false", largeop = "false", movablelimits = "false", accent = "false", lspace = "0.0em", rspace = "0.0em"), Typesetting:-mi("&delta;2", italic = "true", mathvariant = "italic")), Typesetting:-mo("&minus;", mathvariant = "normal", fence = "false", separator = "false", stretchy = "false", symmetric = "false", largeop = "false", movablelimits = "false", accent = "false", lspace = "0.2222222em", rspace = "0.2222222em"), Typesetting:-mrow(Typesetting:-mi("K2", italic = "true", mathvariant = "italic"), Typesetting:-mo("&InvisibleTimes;", mathvariant = "normal", fence = "false", separator = "false", stretchy = "false", symmetric = "false", largeop = "false", movablelimits = "false", accent = "false", lspace = "0.0em", rspace = "0.0em"), Typesetting:-mi("N2", italic = "true", mathvariant = "italic"), Typesetting:-mo("&InvisibleTimes;", mathvariant = "normal", fence = "false", separator = "false", stretchy = "false", symmetric = "false", largeop = "false", movablelimits = "false", accent = "false", lspace = "0.0em", rspace = "0.0em"), Typesetting:-mi("&delta;1", italic = "true", mathvariant = "italic")), Typesetting:-mo("+", mathvariant = "normal", fence = "false", separator = "false", stretchy = "false", symmetric = "false", largeop = "false", movablelimits = "false", accent = "false", lspace = "0.2222222em", rspace = "0.2222222em"), Typesetting:-mrow(Typesetting:-mn("2", mathvariant = "normal"), Typesetting:-mo("&InvisibleTimes;", mathvariant = "normal", fence = "false", separator = "false", stretchy = "false", symmetric = "false", largeop = "false", movablelimits = "false", accent = "false", lspace = "0.0em", rspace = "0.0em"), Typesetting:-mi("K2", italic = "true", mathvariant = "italic"), Typesetting:-mo("&InvisibleTimes;", mathvariant = "normal", fence = "false", separator = "false", stretchy = "false", symmetric = "false", largeop = "false", movablelimits = "false", accent = "false", lspace = "0.0em", rspace = "0.0em"), Typesetting:-mi("N1", italic = "true", mathvariant = "italic"))), mathvariant = "normal")), Typesetting:-mrow(Typesetting:-msup(Typesetting:-mi("K2", italic = "true", mathvariant = "italic"), Typesetting:-mn("2", mathvariant = "normal"), superscriptshift = "0"), Typesetting:-mo("&InvisibleTimes;", mathvariant = "normal", fence = "false", separator = "false", stretchy = "false", symmetric = "false", largeop = "false", movablelimits = "false", accent = "false", lspace = "0.0em", rspace = "0.0em"), Typesetting:-mi("K1", italic = "true", mathvariant = "italic")), linethickness = "1", denomalign = "center", numalign = "center", bevelled = "false")))}{\[\displaystyle B ~\mbox {{\tt }} ^{{\it x2} }~_{{\it x1} }=\frac{1}{4}~\frac{\lambda ^{2}~\delta \mbox {{\tt 2}}\\
\mbox{} ~{\it N2} ~\left({\it K1}\\
\mbox{} ~{\it N1} ~\delta \mbox {{\tt 2}}\\
\mbox{} -{\it K2} ~{\it N2} ~\delta \mbox {{\tt 1}}\\
\mbox{} +2~{\it K2} ~{\it N1} \right)}{{\it K2} ^{2}~{\it K1}\\
\mbox{} }\]}
\end{maplelatex}
\mapleresult
\begin{maplelatex}
\mapleinline{inert}{2d}{Typesetting:-mrow(Typesetting:-mrow(Typesetting:-mi("B", italic = "true", mathvariant = "italic"), Typesetting:-mo("&InvisibleTimes;", mathvariant = "normal", fence = "false", separator = "false", stretchy = "false", symmetric = "false", largeop = "false", movablelimits = "false", accent = "false", lspace = "0.0em", rspace = "0.0em"), Typesetting:-msup(Typesetting:-mi(" ", italic = "true", mathvariant = "italic"), Typesetting:-mi("x2", italic = "true", mathvariant = "italic"), superscriptshift = "0"), Typesetting:-mo("&InvisibleTimes;", mathvariant = "normal", fence = "false", separator = "false", stretchy = "false", symmetric = "false", largeop = "false", movablelimits = "false", accent = "false", lspace = "0.0em", rspace = "0.0em"), Typesetting:-msub(Typesetting:-mi("", italic = "true", mathvariant = "italic"), Typesetting:-mrow(Typesetting:-mi("x2", italic = "true", mathvariant = "italic")), subscriptshift = "0")), Typesetting:-mo("=", mathvariant = "normal", fence = "false", separator = "false", stretchy = "false", symmetric = "false", largeop = "false", movablelimits = "false", accent = "false", lspace = "0.2777778em", rspace = "0.2777778em"), Typesetting:-mrow(Typesetting:-mo("&uminus0;", mathvariant = "normal", fence = "false", separator = "false", stretchy = "false", symmetric = "false", largeop = "false", movablelimits = "false", accent = "false", lspace = "0.2222222em", rspace = "0.2222222em"), Typesetting:-mrow(Typesetting:-mfrac(Typesetting:-mn("1", mathvariant = "normal"), Typesetting:-mn("4", mathvariant = "normal"), linethickness = "1", denomalign = "center", numalign = "center", bevelled = "false"), Typesetting:-mo("&InvisibleTimes;", mathvariant = "normal", fence = "false", separator = "false", stretchy = "false", symmetric = "false", largeop = "false", movablelimits = "false", accent = "false", lspace = "0.0em", rspace = "0.0em"), Typesetting:-mfrac(Typesetting:-mrow(Typesetting:-msup(Typesetting:-mi("&lambda;", italic = "false", mathvariant = "normal"), Typesetting:-mn("2", mathvariant = "normal"), superscriptshift = "0"), Typesetting:-mo("&InvisibleTimes;", mathvariant = "normal", fence = "false", separator = "false", stretchy = "false", symmetric = "false", largeop = "false", movablelimits = "false", accent = "false", lspace = "0.0em", rspace = "0.0em"), Typesetting:-mi("&delta;2", italic = "true", mathvariant = "italic"), Typesetting:-mo("&InvisibleTimes;", mathvariant = "normal", fence = "false", separator = "false", stretchy = "false", symmetric = "false", largeop = "false", movablelimits = "false", accent = "false", lspace = "0.0em", rspace = "0.0em"), Typesetting:-mi("N1", italic = "true", mathvariant = "italic"), Typesetting:-mo("&InvisibleTimes;", mathvariant = "normal", fence = "false", separator = "false", stretchy = "false", symmetric = "false", largeop = "false", movablelimits = "false", accent = "false", lspace = "0.0em", rspace = "0.0em"), Typesetting:-mfenced(Typesetting:-mrow(Typesetting:-mrow(Typesetting:-mi("K1", italic = "true", mathvariant = "italic"), Typesetting:-mo("&InvisibleTimes;", mathvariant = "normal", fence = "false", separator = "false", stretchy = "false", symmetric = "false", largeop = "false", movablelimits = "false", accent = "false", lspace = "0.0em", rspace = "0.0em"), Typesetting:-mi("N1", italic = "true", mathvariant = "italic"), Typesetting:-mo("&InvisibleTimes;", mathvariant = "normal", fence = "false", separator = "false", stretchy = "false", symmetric = "false", largeop = "false", movablelimits = "false", accent = "false", lspace = "0.0em", rspace = "0.0em"), Typesetting:-mi("&delta;2", italic = "true", mathvariant = "italic")), Typesetting:-mo("&minus;", mathvariant = "normal", fence = "false", separator = "false", stretchy = "false", symmetric = "false", largeop = "false", movablelimits = "false", accent = "false", lspace = "0.2222222em", rspace = "0.2222222em"), Typesetting:-mrow(Typesetting:-mi("K2", italic = "true", mathvariant = "italic"), Typesetting:-mo("&InvisibleTimes;", mathvariant = "normal", fence = "false", separator = "false", stretchy = "false", symmetric = "false", largeop = "false", movablelimits = "false", accent = "false", lspace = "0.0em", rspace = "0.0em"), Typesetting:-mi("N2", italic = "true", mathvariant = "italic"), Typesetting:-mo("&InvisibleTimes;", mathvariant = "normal", fence = "false", separator = "false", stretchy = "false", symmetric = "false", largeop = "false", movablelimits = "false", accent = "false", lspace = "0.0em", rspace = "0.0em"), Typesetting:-mi("&delta;1", italic = "true", mathvariant = "italic")), Typesetting:-mo("+", mathvariant = "normal", fence = "false", separator = "false", stretchy = "false", symmetric = "false", largeop = "false", movablelimits = "false", accent = "false", lspace = "0.2222222em", rspace = "0.2222222em"), Typesetting:-mrow(Typesetting:-mn("2", mathvariant = "normal"), Typesetting:-mo("&InvisibleTimes;", mathvariant = "normal", fence = "false", separator = "false", stretchy = "false", symmetric = "false", largeop = "false", movablelimits = "false", accent = "false", lspace = "0.0em", rspace = "0.0em"), Typesetting:-mi("K2", italic = "true", mathvariant = "italic"), Typesetting:-mo("&InvisibleTimes;", mathvariant = "normal", fence = "false", separator = "false", stretchy = "false", symmetric = "false", largeop = "false", movablelimits = "false", accent = "false", lspace = "0.0em", rspace = "0.0em"), Typesetting:-mi("N1", italic = "true", mathvariant = "italic"))), mathvariant = "normal")), Typesetting:-mrow(Typesetting:-msup(Typesetting:-mi("K2", italic = "true", mathvariant = "italic"), Typesetting:-mn("2", mathvariant = "normal"), superscriptshift = "0"), Typesetting:-mo("&InvisibleTimes;", mathvariant = "normal", fence = "false", separator = "false", stretchy = "false", symmetric = "false", largeop = "false", movablelimits = "false", accent = "false", lspace = "0.0em", rspace = "0.0em"), Typesetting:-mi("K1", italic = "true", mathvariant = "italic")), linethickness = "1", denomalign = "center", numalign = "center", bevelled = "false"))))}{\[\displaystyle B ~\mbox {{\tt}} ^{{\it x2} }~_{{\it x2} }=-\frac{1}{4}~\frac{\lambda ^{2}~\delta \mbox {{\tt 2}}\\
\mbox{} ~{\it N1} ~\left({\it K1} ~{\it N1} ~\delta \mbox {{\tt 2}}\\
\mbox{} -{\it K2} ~{\it N2} ~\delta \mbox {{\tt 1}}\\
\mbox{} +2~{\it K2} ~{\it N1} \right)}{{\it K2} ^{2}~{\it K1} }\]}
\end{maplelatex}
\end{maplegroup}
\begin{maplegroup}
\begin{mapleinput}
\mapleinline{active}{2d}{HH := convert(B[i, -j], matrix); 1}{\[\]}
\end{mapleinput}
\begin{maplelatex}
\mapleinline{inert}{2d}{Eingenvalues of the second the second KCC-invariant:}{\[\displaystyle \mbox {{\tt Eingenvalues of the second the second KCC-invariant:}}\]}
\end{maplelatex}
\mapleresult
\begin{maplelatex}
\mapleinline{inert}{2d}{HH := Matrix(id = 18446744074633140214)}{\[\displaystyle \]}
\end{maplelatex}
\end{maplegroup}
\begin{maplegroup}
\begin{mapleinput}
\mapleinline{active}{2d}{eigenvalues(HH); 1}{\[\]}
\end{mapleinput}
\mapleresult
\begin{maplelatex}
\mapleinline{inert}{2d}{0, -\frac{1}{4}*\lambda^2*(K1^2*N1^2*\delta2^2-2*K1*K2*N1*N2*\delta1*\delta2+K2^2*N2^2*\delta1^2+2*K1*K2*N1^2*\delta2+2*K1*K2*N2^2*\delta1)/(K1^2*K2^2)}{\[\displaystyle 0,\,-\frac{1}{4}\,{\frac {{\lambda}^{2} \left( {{\it K1}}^{2}{{\it N1}}^{2}{{\it \delta2}}^{2}-2\,{\it K1}\,{\it K2}\,{\it N1}\,{\it N2}\,{\it \delta1}\,{\it \delta2}\\
\mbox{}+{{\it K2}}^{2}{{\it N2}}^{2}{{\it \delta1}}^{2}+2\,{\it K1}\,{\it K2}\,{{\it N1}}^{2}{\it \delta2}+2\,{\it K1}\,{\it K2}\,{{\it N2}}^{2}{\it \delta1} \right) }{{{\it K1}}^{2}{{\it K2}}^{2}}}\]}
\end{maplelatex}
\end{maplegroup}
\begin{maplegroup}
\begin{mapleinput}
\mapleinline{active}{2d}{show(NLR[i, -j, -k]); 1}{\[\]}
\end{mapleinput}
\begin{maplelatex}
\mapleinline{inert}{2d}{Third KCC-invariant:}{\[\displaystyle \mbox {{\tt Third KCC-invariant:}}\]}
\end{maplelatex}
\mapleresult
\begin{maplelatex}
\mapleinline{inert}{2d}{Typesetting:-mrow(Typesetting:-mrow(Typesetting:-mi("NLR", italic = "true", mathvariant = "italic"), Typesetting:-mo("&InvisibleTimes;", mathvariant = "normal", fence = "false", separator = "false", stretchy = "false", symmetric = "false", largeop = "false", movablelimits = "false", accent = "false", lspace = "0.0em", rspace = "0.0em"), Typesetting:-msup(Typesetting:-mi(" ", italic = "true", mathvariant = "italic"), Typesetting:-mi("x1", italic = "true", mathvariant = "italic"), superscriptshift = "0"), Typesetting:-mo("&InvisibleTimes;", mathvariant = "normal", fence = "false", separator = "false", stretchy = "false", symmetric = "false", largeop = "false", movablelimits = "false", accent = "false", lspace = "0.0em", rspace = "0.0em"), Typesetting:-msub(Typesetting:-mi("", italic = "true", mathvariant = "italic"), Typesetting:-mrow(Typesetting:-mi("x1", italic = "true", mathvariant = "italic")), subscriptshift = "0"), Typesetting:-mo("&InvisibleTimes;", mathvariant = "normal", fence = "false", separator = "false", stretchy = "false", symmetric = "false", largeop = "false", movablelimits = "false", accent = "false", lspace = "0.0em", rspace = "0.0em"), Typesetting:-msub(Typesetting:-mi("", italic = "true", mathvariant = "italic"), Typesetting:-mrow(Typesetting:-mi("x2", italic = "true", mathvariant = "italic")), subscriptshift = "0")), Typesetting:-mo("=", mathvariant = "normal", fence = "false", separator = "false", stretchy = "false", symmetric = "false", largeop = "false", movablelimits = "false", accent = "false", lspace = "0.2777778em", rspace = "0.2777778em"), Typesetting:-mrow(Typesetting:-mfrac(Typesetting:-mn("1", mathvariant = "normal"), Typesetting:-mn("4", mathvariant = "normal"), linethickness = "1", denomalign = "center", numalign = "center", bevelled = "false"), Typesetting:-mo("&InvisibleTimes;", mathvariant = "normal", fence = "false", separator = "false", stretchy = "false", symmetric = "false", largeop = "false", movablelimits = "false", accent = "false", lspace = "0.0em", rspace = "0.0em"), Typesetting:-mfrac(Typesetting:-mrow(Typesetting:-msup(Typesetting:-mi("&lambda;", italic = "false", mathvariant = "normal"), Typesetting:-mn("2", mathvariant = "normal"), superscriptshift = "0"), Typesetting:-mo("&InvisibleTimes;", mathvariant = "normal", fence = "false", separator = "false", stretchy = "false", symmetric = "false", largeop = "false", movablelimits = "false", accent = "false", lspace = "0.0em", rspace = "0.0em"), Typesetting:-mi("&delta;1", italic = "true", mathvariant = "italic"), Typesetting:-mo("&InvisibleTimes;", mathvariant = "normal", fence = "false", separator = "false", stretchy = "false", symmetric = "false", largeop = "false", movablelimits = "false", accent = "false", lspace = "0.0em", rspace = "0.0em"), Typesetting:-mfenced(Typesetting:-mrow(Typesetting:-mrow(Typesetting:-mi("K1", italic = "true", mathvariant = "italic"), Typesetting:-mo("&InvisibleTimes;", mathvariant = "normal", fence = "false", separator = "false", stretchy = "false", symmetric = "false", largeop = "false", movablelimits = "false", accent = "false", lspace = "0.0em", rspace = "0.0em"), Typesetting:-mi("N1", italic = "true", mathvariant = "italic"), Typesetting:-mo("&InvisibleTimes;", mathvariant = "normal", fence = "false", separator = "false", stretchy = "false", symmetric = "false", largeop = "false", movablelimits = "false", accent = "false", lspace = "0.0em", rspace = "0.0em"), Typesetting:-mi("&delta;2", italic = "true", mathvariant = "italic")), Typesetting:-mo("&minus;", mathvariant = "normal", fence = "false", separator = "false", stretchy = "false", symmetric = "false", largeop = "false", movablelimits = "false", accent = "false", lspace = "0.2222222em", rspace = "0.2222222em"), Typesetting:-mrow(Typesetting:-mi("K2", italic = "true", mathvariant = "italic"), Typesetting:-mo("&InvisibleTimes;", mathvariant = "normal", fence = "false", separator = "false", stretchy = "false", symmetric = "false", largeop = "false", movablelimits = "false", accent = "false", lspace = "0.0em", rspace = "0.0em"), Typesetting:-mi("N2", italic = "true", mathvariant = "italic"), Typesetting:-mo("&InvisibleTimes;", mathvariant = "normal", fence = "false", separator = "false", stretchy = "false", symmetric = "false", largeop = "false", movablelimits = "false", accent = "false", lspace = "0.0em", rspace = "0.0em"), Typesetting:-mi("&delta;1", italic = "true", mathvariant = "italic")), Typesetting:-mo("&minus;", mathvariant = "normal", fence = "false", separator = "false", stretchy = "false", symmetric = "false", largeop = "false", movablelimits = "false", accent = "false", lspace = "0.2222222em", rspace = "0.2222222em"), Typesetting:-mrow(Typesetting:-mn("2", mathvariant = "normal"), Typesetting:-mo("&InvisibleTimes;", mathvariant = "normal", fence = "false", separator = "false", stretchy = "false", symmetric = "false", largeop = "false", movablelimits = "false", accent = "false", lspace = "0.0em", rspace = "0.0em"), Typesetting:-mi("K1", italic = "true", mathvariant = "italic"), Typesetting:-mo("&InvisibleTimes;", mathvariant = "normal", fence = "false", separator = "false", stretchy = "false", symmetric = "false", largeop = "false", movablelimits = "false", accent = "false", lspace = "0.0em", rspace = "0.0em"), Typesetting:-mi("N2", italic = "true", mathvariant = "italic"))), mathvariant = "normal")), Typesetting:-mrow(Typesetting:-msup(Typesetting:-mi("K1", italic = "true", mathvariant = "italic"), Typesetting:-mn("2", mathvariant = "normal"), superscriptshift = "0"), Typesetting:-mo("&InvisibleTimes;", mathvariant = "normal", fence = "false", separator = "false", stretchy = "false", symmetric = "false", largeop = "false", movablelimits = "false", accent = "false", lspace = "0.0em", rspace = "0.0em"), Typesetting:-mi("K2", italic = "true", mathvariant = "italic")), linethickness = "1", denomalign = "center", numalign = "center", bevelled = "false")))}
{\[\displaystyle {\it NLR} ~\mbox {{\tt }} ^{{\it x1} }~_{{\it x1} }~_{{\it x2} }=\frac{1}{4}~\frac{\lambda ^{2}~\delta \mbox {{\tt 1}}\\
\mbox{} ~\left({\it K1} ~{\it N1}\\
\mbox{} ~\delta \mbox {{\tt 2}} -{\it K2}\\
\mbox{} ~{\it N2} ~\delta \mbox {{\tt 1}}\\
\mbox{} -2~{\it K1} ~{\it N2} \right)}{{\it K1} ^{2}~{\it K2}\\
\mbox{} }\]}
\end{maplelatex}
\mapleresult

\begin{maplelatex}
\mapleinline{inert}{2d}{Typesetting:-mrow(Typesetting:-mrow(Typesetting:-mi("NLR", italic = "true", mathvariant = "italic"), Typesetting:-mo("&InvisibleTimes;", mathvariant = "normal", fence = "false", separator = "false", stretchy = "false", symmetric = "false", largeop = "false", movablelimits = "false", accent = "false", lspace = "0.0em", rspace = "0.0em"), Typesetting:-msup(Typesetting:-mi(" ", italic = "true", mathvariant = "italic"), Typesetting:-mi("x2", italic = "true", mathvariant = "italic"), superscriptshift = "0"), Typesetting:-mo("&InvisibleTimes;", mathvariant = "normal", fence = "false", separator = "false", stretchy = "false", symmetric = "false", largeop = "false", movablelimits = "false", accent = "false", lspace = "0.0em", rspace = "0.0em"), Typesetting:-msub(Typesetting:-mi("", italic = "true", mathvariant = "italic"), Typesetting:-mrow(Typesetting:-mi("x1", italic = "true", mathvariant = "italic")), subscriptshift = "0"), Typesetting:-mo("&InvisibleTimes;", mathvariant = "normal", fence = "false", separator = "false", stretchy = "false", symmetric = "false", largeop = "false", movablelimits = "false", accent = "false", lspace = "0.0em", rspace = "0.0em"), Typesetting:-msub(Typesetting:-mi("", italic = "true", mathvariant = "italic"), Typesetting:-mrow(Typesetting:-mi("x2", italic = "true", mathvariant = "italic")), subscriptshift = "0")), Typesetting:-mo("=", mathvariant = "normal", fence = "false", separator = "false", stretchy = "false", symmetric = "false", largeop = "false", movablelimits = "false", accent = "false", lspace = "0.2777778em", rspace = "0.2777778em"), Typesetting:-mrow(Typesetting:-mfrac(Typesetting:-mn("1", mathvariant = "normal"), Typesetting:-mn("4", mathvariant = "normal"), linethickness = "1", denomalign = "center", numalign = "center", bevelled = "false"), Typesetting:-mo("&InvisibleTimes;", mathvariant = "normal", fence = "false", separator = "false", stretchy = "false", symmetric = "false", largeop = "false", movablelimits = "false", accent = "false", lspace = "0.0em", rspace = "0.0em"), Typesetting:-mfrac(Typesetting:-mrow(Typesetting:-msup(Typesetting:-mi("&lambda;", italic = "false", mathvariant = "normal"), Typesetting:-mn("2", mathvariant = "normal"), superscriptshift = "0"), Typesetting:-mo("&InvisibleTimes;", mathvariant = "normal", fence = "false", separator = "false", stretchy = "false", symmetric = "false", largeop = "false", movablelimits = "false", accent = "false", lspace = "0.0em", rspace = "0.0em"), Typesetting:-mi("&delta;2", italic = "true", mathvariant = "italic"), Typesetting:-mo("&InvisibleTimes;", mathvariant = "normal", fence = "false", separator = "false", stretchy = "false", symmetric = "false", largeop = "false", movablelimits = "false", accent = "false", lspace = "0.0em", rspace = "0.0em"), Typesetting:-mfenced(Typesetting:-mrow(Typesetting:-mrow(Typesetting:-mi("K1", italic = "true", mathvariant = "italic"), Typesetting:-mo("&InvisibleTimes;", mathvariant = "normal", fence = "false", separator = "false", stretchy = "false", symmetric = "false", largeop = "false", movablelimits = "false", accent = "false", lspace = "0.0em", rspace = "0.0em"), Typesetting:-mi("N1", italic = "true", mathvariant = "italic"), Typesetting:-mo("&InvisibleTimes;", mathvariant = "normal", fence = "false", separator = "false", stretchy = "false", symmetric = "false", largeop = "false", movablelimits = "false", accent = "false", lspace = "0.0em", rspace = "0.0em"), Typesetting:-mi("&delta;2", italic = "true", mathvariant = "italic")), Typesetting:-mo("&minus;", mathvariant = "normal", fence = "false", separator = "false", stretchy = "false", symmetric = "false", largeop = "false", movablelimits = "false", accent = "false", lspace = "0.2222222em", rspace = "0.2222222em"), Typesetting:-mrow(Typesetting:-mi("K2", italic = "true", mathvariant = "italic"), Typesetting:-mo("&InvisibleTimes;", mathvariant = "normal", fence = "false", separator = "false", stretchy = "false", symmetric = "false", largeop = "false", movablelimits = "false", accent = "false", lspace = "0.0em", rspace = "0.0em"), Typesetting:-mi("N2", italic = "true", mathvariant = "italic"), Typesetting:-mo("&InvisibleTimes;", mathvariant = "normal", fence = "false", separator = "false", stretchy = "false", symmetric = "false", largeop = "false", movablelimits = "false", accent = "false", lspace = "0.0em", rspace = "0.0em"), Typesetting:-mi("&delta;1", italic = "true", mathvariant = "italic")), Typesetting:-mo("+", mathvariant = "normal", fence = "false", separator = "false", stretchy = "false", symmetric = "false", largeop = "false", movablelimits = "false", accent = "false", lspace = "0.2222222em", rspace = "0.2222222em"), Typesetting:-mrow(Typesetting:-mn("2", mathvariant = "normal"), Typesetting:-mo("&InvisibleTimes;", mathvariant = "normal", fence = "false", separator = "false", stretchy = "false", symmetric = "false", largeop = "false", movablelimits = "false", accent = "false", lspace = "0.0em", rspace = "0.0em"), Typesetting:-mi("K2", italic = "true", mathvariant = "italic"), Typesetting:-mo("&InvisibleTimes;", mathvariant = "normal", fence = "false", separator = "false", stretchy = "false", symmetric = "false", largeop = "false", movablelimits = "false", accent = "false", lspace = "0.0em", rspace = "0.0em"), Typesetting:-mi("N1", italic = "true", mathvariant = "italic"))), mathvariant = "normal")), Typesetting:-mrow(Typesetting:-msup(Typesetting:-mi("K2", italic = "true", mathvariant = "italic"), Typesetting:-mn("2", mathvariant = "normal"), superscriptshift = "0"), Typesetting:-mo("&InvisibleTimes;", mathvariant = "normal", fence = "false", separator = "false", stretchy = "false", symmetric = "false", largeop = "false", movablelimits = "false", accent = "false", lspace = "0.0em", rspace = "0.0em"), Typesetting:-mi("K1", italic = "true", mathvariant = "italic")), linethickness = "1", denomalign = "center", numalign = "center", bevelled = "false")))}{\[\displaystyle {\it NLR} ~\mbox {{\tt }} ^{{\it x2} }~_{{\it x1} }~_{{\it x2} }=\frac{1}{4}~\frac{\lambda ^{2}~\delta \mbox {{\tt 2}}\\
\mbox{} ~\left({\it K1} ~{\it N1}\\
\mbox{} ~\delta \mbox {{\tt 2}}\\
\mbox{} -{\it K2} ~{\it N2} ~\delta \mbox {{\tt 1}}\\
\mbox{} +2~{\it K2} ~{\it N1}\\
\mbox{} \right)}{{\it K2} ^{2}~{\it K1} }\]}
\end{maplelatex}
\end{maplegroup}
\begin{maplegroup}
\begin{mapleinput}
\mapleinline{active}{2d}{show(B[i, -j, -k, -l]); 1}{\[\]}
\end{mapleinput}
\begin{maplelatex}
\mapleinline{inert}{2d}{Fourth KCC-invariant:}{\[\displaystyle \mbox {{\tt Fourth KCC-invariant:}}\]}
\end{maplelatex}
\mapleresult
\begin{maplelatex}
\mapleinline{inert}{2d}{Typesetting:-mrow(Typesetting:-mrow(Typesetting:-mi("B", italic = "true", mathvariant = "italic"), Typesetting:-mo("&InvisibleTimes;", mathvariant = "normal", fence = "false", separator = "false", stretchy = "false", symmetric = "false", largeop = "false", movablelimits = "false", accent = "false", lspace = "0.0em", rspace = "0.0em"), Typesetting:-msup(Typesetting:-mi(" ", italic = "true", mathvariant = "italic"), Typesetting:-mi("x1", italic = "true", mathvariant = "italic"), superscriptshift = "0"), Typesetting:-mo("&InvisibleTimes;", mathvariant = "normal", fence = "false", separator = "false", stretchy = "false", symmetric = "false", largeop = "false", movablelimits = "false", accent = "false", lspace = "0.0em", rspace = "0.0em"), Typesetting:-msub(Typesetting:-mi("", italic = "true", mathvariant = "italic"), Typesetting:-mrow(Typesetting:-mi("x1", italic = "true", mathvariant = "italic")), subscriptshift = "0"), Typesetting:-mo("&InvisibleTimes;", mathvariant = "normal", fence = "false", separator = "false", stretchy = "false", symmetric = "false", largeop = "false", movablelimits = "false", accent = "false", lspace = "0.0em", rspace = "0.0em"), Typesetting:-msub(Typesetting:-mi("", italic = "true", mathvariant = "italic"), Typesetting:-mrow(Typesetting:-mi("x1", italic = "true", mathvariant = "italic")), subscriptshift = "0"), Typesetting:-mo("&InvisibleTimes;", mathvariant = "normal", fence = "false", separator = "false", stretchy = "false", symmetric = "false", largeop = "false", movablelimits = "false", accent = "false", lspace = "0.0em", rspace = "0.0em"), Typesetting:-msub(Typesetting:-mi("", italic = "true", mathvariant = "italic"), Typesetting:-mrow(Typesetting:-mi("x2", italic = "true", mathvariant = "italic")), subscriptshift = "0")), Typesetting:-mo("=", mathvariant = "normal", fence = "false", separator = "false", stretchy = "false", symmetric = "false", largeop = "false", movablelimits = "false", accent = "false", lspace = "0.2777778em", rspace = "0.2777778em"), Typesetting:-mrow(Typesetting:-mfrac(Typesetting:-mn("1", mathvariant = "normal"), Typesetting:-mn("4", mathvariant = "normal"), linethickness = "1", denomalign = "center", numalign = "center", bevelled = "false"), Typesetting:-mo("&InvisibleTimes;", mathvariant = "normal", fence = "false", separator = "false", stretchy = "false", symmetric = "false", largeop = "false", movablelimits = "false", accent = "false", lspace = "0.0em", rspace = "0.0em"), Typesetting:-mfrac(Typesetting:-mrow(Typesetting:-msup(Typesetting:-mi("&lambda;", italic = "false", mathvariant = "normal"), Typesetting:-mn("2", mathvariant = "normal"), superscriptshift = "0"), Typesetting:-mo("&InvisibleTimes;", mathvariant = "normal", fence = "false", separator = "false", stretchy = "false", symmetric = "false", largeop = "false", movablelimits = "false", accent = "false", lspace = "0.0em", rspace = "0.0em"), Typesetting:-mi("&delta;1", italic = "true", mathvariant = "italic"), Typesetting:-mo("&InvisibleTimes;", mathvariant = "normal", fence = "false", separator = "false", stretchy = "false", symmetric = "false", largeop = "false", movablelimits = "false", accent = "false", lspace = "0.0em", rspace = "0.0em"), Typesetting:-mi("&delta;2", italic = "true", mathvariant = "italic")), Typesetting:-mrow(Typesetting:-mi("K1", italic = "true", mathvariant = "italic"), Typesetting:-mo("&InvisibleTimes;", mathvariant = "normal", fence = "false", separator = "false", stretchy = "false", symmetric = "false", largeop = "false", movablelimits = "false", accent = "false", lspace = "0.0em", rspace = "0.0em"), Typesetting:-mi("K2", italic = "true", mathvariant = "italic")), linethickness = "1", denomalign = "center", numalign = "center", bevelled = "false")))}
{\[\displaystyle B ~\mbox {{\tt }} ^{{\it x1} }~_{{\it x1} }~_{{\it x1} }~_{{\it x2} }=\frac{1}{4}~\frac{\lambda ^{2}~\delta \mbox {{\tt 1}}\\
\mbox{} ~\delta \mbox {{\tt 2}}\\
\mbox{} }{{\it K1} ~{\it K2}\\
\mbox{} }\]}
\end{maplelatex}
\mapleresult
\begin{maplelatex}
\mapleinline{inert}{2d}{Typesetting:-mrow(Typesetting:-mrow(Typesetting:-mi("B", italic = "true", mathvariant = "italic"), Typesetting:-mo("&InvisibleTimes;", mathvariant = "normal", fence = "false", separator = "false", stretchy = "false", symmetric = "false", largeop = "false", movablelimits = "false", accent = "false", lspace = "0.0em", rspace = "0.0em"), Typesetting:-msup(Typesetting:-mi(" ", italic = "true", mathvariant = "italic"), Typesetting:-mi("x1", italic = "true", mathvariant = "italic"), superscriptshift = "0"), Typesetting:-mo("&InvisibleTimes;", mathvariant = "normal", fence = "false", separator = "false", stretchy = "false", symmetric = "false", largeop = "false", movablelimits = "false", accent = "false", lspace = "0.0em", rspace = "0.0em"), Typesetting:-msub(Typesetting:-mi("", italic = "true", mathvariant = "italic"), Typesetting:-mrow(Typesetting:-mi("x2", italic = "true", mathvariant = "italic")), subscriptshift = "0"), Typesetting:-mo("&InvisibleTimes;", mathvariant = "normal", fence = "false", separator = "false", stretchy = "false", symmetric = "false", largeop = "false", movablelimits = "false", accent = "false", lspace = "0.0em", rspace = "0.0em"), Typesetting:-msub(Typesetting:-mi("", italic = "true", mathvariant = "italic"), Typesetting:-mrow(Typesetting:-mi("x1", italic = "true", mathvariant = "italic")), subscriptshift = "0"), Typesetting:-mo("&InvisibleTimes;", mathvariant = "normal", fence = "false", separator = "false", stretchy = "false", symmetric = "false", largeop = "false", movablelimits = "false", accent = "false", lspace = "0.0em", rspace = "0.0em"), Typesetting:-msub(Typesetting:-mi("", italic = "true", mathvariant = "italic"), Typesetting:-mrow(Typesetting:-mi("x2", italic = "true", mathvariant = "italic")), subscriptshift = "0")), Typesetting:-mo("=", mathvariant = "normal", fence = "false", separator = "false", stretchy = "false", symmetric = "false", largeop = "false", movablelimits = "false", accent = "false", lspace = "0.2777778em", rspace = "0.2777778em"), Typesetting:-mrow(Typesetting:-mo("&uminus0;", mathvariant = "normal", fence = "false", separator = "false", stretchy = "false", symmetric = "false", largeop = "false", movablelimits = "false", accent = "false", lspace = "0.2222222em", rspace = "0.2222222em"), Typesetting:-mrow(Typesetting:-mfrac(Typesetting:-mn("1", mathvariant = "normal"), Typesetting:-mn("4", mathvariant = "normal"), linethickness = "1", denomalign = "center", numalign = "center", bevelled = "false"), Typesetting:-mo("&InvisibleTimes;", mathvariant = "normal", fence = "false", separator = "false", stretchy = "false", symmetric = "false", largeop = "false", movablelimits = "false", accent = "false", lspace = "0.0em", rspace = "0.0em"), Typesetting:-mfrac(Typesetting:-mrow(Typesetting:-msup(Typesetting:-mi("&lambda;", italic = "false", mathvariant = "normal"), Typesetting:-mn("2", mathvariant = "normal"), superscriptshift = "0"), Typesetting:-mo("&InvisibleTimes;", mathvariant = "normal", fence = "false", separator = "false", stretchy = "false", symmetric = "false", largeop = "false", movablelimits = "false", accent = "false", lspace = "0.0em", rspace = "0.0em"), Typesetting:-mi("&delta;1", italic = "true", mathvariant = "italic"), Typesetting:-mo("&InvisibleTimes;", mathvariant = "normal", fence = "false", separator = "false", stretchy = "false", symmetric = "false", largeop = "false", movablelimits = "false", accent = "false", lspace = "0.0em", rspace = "0.0em"), Typesetting:-mfenced(Typesetting:-mrow(Typesetting:-mrow(Typesetting:-mi("K2", italic = "true", mathvariant = "italic"), Typesetting:-mo("&InvisibleTimes;", mathvariant = "normal", fence = "false", separator = "false", stretchy = "false", symmetric = "false", largeop = "false", movablelimits = "false", accent = "false", lspace = "0.0em", rspace = "0.0em"), Typesetting:-mi("&delta;1", italic = "true", mathvariant = "italic")), Typesetting:-mo("+", mathvariant = "normal", fence = "false", separator = "false", stretchy = "false", symmetric = "false", largeop = "false", movablelimits = "false", accent = "false", lspace = "0.2222222em", rspace = "0.2222222em"), Typesetting:-mrow(Typesetting:-mn("2", mathvariant = "normal"), Typesetting:-mo("&InvisibleTimes;", mathvariant = "normal", fence = "false", separator = "false", stretchy = "false", symmetric = "false", largeop = "false", movablelimits = "false", accent = "false", lspace = "0.0em", rspace = "0.0em"), Typesetting:-mi("K1", italic = "true", mathvariant = "italic"))), mathvariant = "normal")), Typesetting:-mrow(Typesetting:-msup(Typesetting:-mi("K1", italic = "true", mathvariant = "italic"), Typesetting:-mn("2", mathvariant = "normal"), superscriptshift = "0"), Typesetting:-mo("&InvisibleTimes;", mathvariant = "normal", fence = "false", separator = "false", stretchy = "false", symmetric = "false", largeop = "false", movablelimits = "false", accent = "false", lspace = "0.0em", rspace = "0.0em"), Typesetting:-mi("K2", italic = "true", mathvariant = "italic")), linethickness = "1", denomalign = "center", numalign = "center", bevelled = "false"))))}{\[\displaystyle B ~\mbox {{\tt }} ^{{\it x1} }~_{{\it x2} }~_{{\it x1} }~_{{\it x2} }=-\frac{1}{4}~\frac{\lambda ^{2}~\delta \mbox {{\tt 1}}\\
\mbox{} ~\left({\it K2} ~\delta \mbox {{\tt 1}}\\
\mbox{} +2~{\it K1}\\
\mbox{} \right)}{{\it K1}\\
\mbox{} ^{2}~{\it K2} }\]}
\end{maplelatex}
\mapleresult
\begin{maplelatex}
\mapleinline{inert}{2d}{Typesetting:-mrow(Typesetting:-mrow(Typesetting:-mi("B", italic = "true", mathvariant = "italic"), Typesetting:-mo("&InvisibleTimes;", mathvariant = "normal", fence = "false", separator = "false", stretchy = "false", symmetric = "false", largeop = "false", movablelimits = "false", accent = "false", lspace = "0.0em", rspace = "0.0em"), Typesetting:-msup(Typesetting:-mi(" ", italic = "true", mathvariant = "italic"), Typesetting:-mi("x2", italic = "true", mathvariant = "italic"), superscriptshift = "0"), Typesetting:-mo("&InvisibleTimes;", mathvariant = "normal", fence = "false", separator = "false", stretchy = "false", symmetric = "false", largeop = "false", movablelimits = "false", accent = "false", lspace = "0.0em", rspace = "0.0em"), Typesetting:-msub(Typesetting:-mi("", italic = "true", mathvariant = "italic"), Typesetting:-mrow(Typesetting:-mi("x1", italic = "true", mathvariant = "italic")), subscriptshift = "0"), Typesetting:-mo("&InvisibleTimes;", mathvariant = "normal", fence = "false", separator = "false", stretchy = "false", symmetric = "false", largeop = "false", movablelimits = "false", accent = "false", lspace = "0.0em", rspace = "0.0em"), Typesetting:-msub(Typesetting:-mi("", italic = "true", mathvariant = "italic"), Typesetting:-mrow(Typesetting:-mi("x1", italic = "true", mathvariant = "italic")), subscriptshift = "0"), Typesetting:-mo("&InvisibleTimes;", mathvariant = "normal", fence = "false", separator = "false", stretchy = "false", symmetric = "false", largeop = "false", movablelimits = "false", accent = "false", lspace = "0.0em", rspace = "0.0em"), Typesetting:-msub(Typesetting:-mi("", italic = "true", mathvariant = "italic"), Typesetting:-mrow(Typesetting:-mi("x2", italic = "true", mathvariant = "italic")), subscriptshift = "0")), Typesetting:-mo("=", mathvariant = "normal", fence = "false", separator = "false", stretchy = "false", symmetric = "false", largeop = "false", movablelimits = "false", accent = "false", lspace = "0.2777778em", rspace = "0.2777778em"), Typesetting:-mrow(Typesetting:-mfrac(Typesetting:-mn("1", mathvariant = "normal"), Typesetting:-mn("4", mathvariant = "normal"), linethickness = "1", denomalign = "center", numalign = "center", bevelled = "false"), Typesetting:-mo("&InvisibleTimes;", mathvariant = "normal", fence = "false", separator = "false", stretchy = "false", symmetric = "false", largeop = "false", movablelimits = "false", accent = "false", lspace = "0.0em", rspace = "0.0em"), Typesetting:-mfrac(Typesetting:-mrow(Typesetting:-msup(Typesetting:-mi("&lambda;", italic = "false", mathvariant = "normal"), Typesetting:-mn("2", mathvariant = "normal"), superscriptshift = "0"), Typesetting:-mo("&InvisibleTimes;", mathvariant = "normal", fence = "false", separator = "false", stretchy = "false", symmetric = "false", largeop = "false", movablelimits = "false", accent = "false", lspace = "0.0em", rspace = "0.0em"), Typesetting:-mi("&delta;2", italic = "true", mathvariant = "italic"), Typesetting:-mo("&InvisibleTimes;", mathvariant = "normal", fence = "false", separator = "false", stretchy = "false", symmetric = "false", largeop = "false", movablelimits = "false", accent = "false", lspace = "0.0em", rspace = "0.0em"), Typesetting:-mfenced(Typesetting:-mrow(Typesetting:-mrow(Typesetting:-mi("K1", italic = "true", mathvariant = "italic"), Typesetting:-mo("&InvisibleTimes;", mathvariant = "normal", fence = "false", separator = "false", stretchy = "false", symmetric = "false", largeop = "false", movablelimits = "false", accent = "false", lspace = "0.0em", rspace = "0.0em"), Typesetting:-mi("&delta;2", italic = "true", mathvariant = "italic")), Typesetting:-mo("+", mathvariant = "normal", fence = "false", separator = "false", stretchy = "false", symmetric = "false", largeop = "false", movablelimits = "false", accent = "false", lspace = "0.2222222em", rspace = "0.2222222em"), Typesetting:-mrow(Typesetting:-mn("2", mathvariant = "normal"), Typesetting:-mo("&InvisibleTimes;", mathvariant = "normal", fence = "false", separator = "false", stretchy = "false", symmetric = "false", largeop = "false", movablelimits = "false", accent = "false", lspace = "0.0em", rspace = "0.0em"), Typesetting:-mi("K2", italic = "true", mathvariant = "italic"))), mathvariant = "normal")), Typesetting:-mrow(Typesetting:-msup(Typesetting:-mi("K2", italic = "true", mathvariant = "italic"), Typesetting:-mn("2", mathvariant = "normal"), superscriptshift = "0"), Typesetting:-mo("&InvisibleTimes;", mathvariant = "normal", fence = "false", separator = "false", stretchy = "false", symmetric = "false", largeop = "false", movablelimits = "false", accent = "false", lspace = "0.0em", rspace = "0.0em"), Typesetting:-mi("K1", italic = "true", mathvariant = "italic")), linethickness = "1", denomalign = "center", numalign = "center", bevelled = "false")))}{\[\displaystyle B ~\mbox {{\tt }} ^{{\it x2} }~_{{\it x1} }~_{{\it x1} }~_{{\it x2} }=\frac{1}{4}~\frac{\lambda ^{2}~\delta \mbox {{\tt 2}}\\
\mbox{} ~\left({\it K1} ~\delta \mbox {{\tt 2}}\\
\mbox{} +2~{\it K2}\\
\mbox{} \right)}{{\it K2}\\
\mbox{} ^{2}~{\it K1} }\]}
\end{maplelatex}
\mapleresult
\begin{maplelatex}
\mapleinline{inert}{2d}{Typesetting:-mrow(Typesetting:-mrow(Typesetting:-mi("B", italic = "true", mathvariant = "italic"), Typesetting:-mo("&InvisibleTimes;", mathvariant = "normal", fence = "false", separator = "false", stretchy = "false", symmetric = "false", largeop = "false", movablelimits = "false", accent = "false", lspace = "0.0em", rspace = "0.0em"), Typesetting:-msup(Typesetting:-mi(" ", italic = "true", mathvariant = "italic"), Typesetting:-mi("x2", italic = "true", mathvariant = "italic"), superscriptshift = "0"), Typesetting:-mo("&InvisibleTimes;", mathvariant = "normal", fence = "false", separator = "false", stretchy = "false", symmetric = "false", largeop = "false", movablelimits = "false", accent = "false", lspace = "0.0em", rspace = "0.0em"), Typesetting:-msub(Typesetting:-mi("", italic = "true", mathvariant = "italic"), Typesetting:-mrow(Typesetting:-mi("x2", italic = "true", mathvariant = "italic")), subscriptshift = "0"), Typesetting:-mo("&InvisibleTimes;", mathvariant = "normal", fence = "false", separator = "false", stretchy = "false", symmetric = "false", largeop = "false", movablelimits = "false", accent = "false", lspace = "0.0em", rspace = "0.0em"), Typesetting:-msub(Typesetting:-mi("", italic = "true", mathvariant = "italic"), Typesetting:-mrow(Typesetting:-mi("x1", italic = "true", mathvariant = "italic")), subscriptshift = "0"), Typesetting:-mo("&InvisibleTimes;", mathvariant = "normal", fence = "false", separator = "false", stretchy = "false", symmetric = "false", largeop = "false", movablelimits = "false", accent = "false", lspace = "0.0em", rspace = "0.0em"), Typesetting:-msub(Typesetting:-mi("", italic = "true", mathvariant = "italic"), Typesetting:-mrow(Typesetting:-mi("x2", italic = "true", mathvariant = "italic")), subscriptshift = "0")), Typesetting:-mo("=", mathvariant = "normal", fence = "false", separator = "false", stretchy = "false", symmetric = "false", largeop = "false", movablelimits = "false", accent = "false", lspace = "0.2777778em", rspace = "0.2777778em"), Typesetting:-mrow(Typesetting:-mo("&uminus0;", mathvariant = "normal", fence = "false", separator = "false", stretchy = "false", symmetric = "false", largeop = "false", movablelimits = "false", accent = "false", lspace = "0.2222222em", rspace = "0.2222222em"), Typesetting:-mrow(Typesetting:-mfrac(Typesetting:-mn("1", mathvariant = "normal"), Typesetting:-mn("4", mathvariant = "normal"), linethickness = "1", denomalign = "center", numalign = "center", bevelled = "false"), Typesetting:-mo("&InvisibleTimes;", mathvariant = "normal", fence = "false", separator = "false", stretchy = "false", symmetric = "false", largeop = "false", movablelimits = "false", accent = "false", lspace = "0.0em", rspace = "0.0em"), Typesetting:-mfrac(Typesetting:-mrow(Typesetting:-msup(Typesetting:-mi("&lambda;", italic = "false", mathvariant = "normal"), Typesetting:-mn("2", mathvariant = "normal"), superscriptshift = "0"), Typesetting:-mo("&InvisibleTimes;", mathvariant = "normal", fence = "false", separator = "false", stretchy = "false", symmetric = "false", largeop = "false", movablelimits = "false", accent = "false", lspace = "0.0em", rspace = "0.0em"), Typesetting:-mi("&delta;1", italic = "true", mathvariant = "italic"), Typesetting:-mo("&InvisibleTimes;", mathvariant = "normal", fence = "false", separator = "false", stretchy = "false", symmetric = "false", largeop = "false", movablelimits = "false", accent = "false", lspace = "0.0em", rspace = "0.0em"), Typesetting:-mi("&delta;2", italic = "true", mathvariant = "italic")), Typesetting:-mrow(Typesetting:-mi("K1", italic = "true", mathvariant = "italic"), Typesetting:-mo("&InvisibleTimes;", mathvariant = "normal", fence = "false", separator = "false", stretchy = "false", symmetric = "false", largeop = "false", movablelimits = "false", accent = "false", lspace = "0.0em", rspace = "0.0em"), Typesetting:-mi("K2", italic = "true", mathvariant = "italic")), linethickness = "1", denomalign = "center", numalign = "center", bevelled = "false"))))}{\[\displaystyle B ~\mbox {{\tt }} ^{{\it x2} }~_{{\it x2} }~_{{\it x1} }~_{{\it x2} }=-\frac{1}{4}~\frac{\lambda ^{2}~\delta \mbox {{\tt 1}}\\
\mbox{} ~\delta \mbox {{\tt 2}}\\
\mbox{} }{{\it K1} ~{\it K2}\\
\mbox{} }\]}
\end{maplelatex}
\end{maplegroup}
\begin{maplegroup}
\begin{mapleinput}
\mapleinline{active}{2d}{show(G[i, -j, -k, -l]); 1}{\[\]}
\end{mapleinput}
\begin{maplelatex}
\mapleinline{inert}{2d}{Fifth KCC-invariant:}{\[\displaystyle \mbox {{\tt Fifth KCC-invariant:}}\]}
\end{maplelatex}
\mapleresult
\begin{maplelatex}
\mapleinline{inert}{2d}{Typesetting:-mrow(Typesetting:-mrow(Typesetting:-mi("G", italic = "true", mathvariant = "italic"), Typesetting:-mo("&InvisibleTimes;", mathvariant = "normal", fence = "false", separator = "false", stretchy = "false", symmetric = "false", largeop = "false", movablelimits = "false", accent = "false", lspace = "0.0em", rspace = "0.0em"), Typesetting:-msup(Typesetting:-mi(" ", italic = "true", mathvariant = "italic"), Typesetting:-mi("i", italic = "true", mathvariant = "italic"), superscriptshift = "0"), Typesetting:-mo("&InvisibleTimes;", mathvariant = "normal", fence = "false", separator = "false", stretchy = "false", symmetric = "false", largeop = "false", movablelimits = "false", accent = "false", lspace = "0.0em", rspace = "0.0em"), Typesetting:-msub(Typesetting:-mi("", italic = "true", mathvariant = "italic"), Typesetting:-mrow(Typesetting:-mi("j", italic = "true", mathvariant = "italic")), subscriptshift = "0"), Typesetting:-mo("&InvisibleTimes;", mathvariant = "normal", fence = "false", separator = "false", stretchy = "false", symmetric = "false", largeop = "false", movablelimits = "false", accent = "false", lspace = "0.0em", rspace = "0.0em"), Typesetting:-msub(Typesetting:-mi("", italic = "true", mathvariant = "italic"), Typesetting:-mrow(Typesetting:-mi("k", italic = "true", mathvariant = "italic")), subscriptshift = "0"), Typesetting:-mo("&InvisibleTimes;", mathvariant = "normal", fence = "false", separator = "false", stretchy = "false", symmetric = "false", largeop = "false", movablelimits = "false", accent = "false", lspace = "0.0em", rspace = "0.0em"), Typesetting:-msub(Typesetting:-mi("", italic = "true", mathvariant = "italic"), Typesetting:-mrow(Typesetting:-mi("l", italic = "true", mathvariant = "italic")), subscriptshift = "0")), Typesetting:-mo("=", mathvariant = "normal", fence = "false", separator = "false", stretchy = "false", symmetric = "false", largeop = "false", movablelimits = "false", accent = "false", lspace = "0.2777778em", rspace = "0.2777778em"), Typesetting:-mn("0", mathvariant = "normal"))}{\[\displaystyle G ~\mbox {{\tt }} ^{i }~_{j }~_{k }~_{l }=0\]}
\end{maplelatex}
\end{maplegroup}

\mapleinline{inert}{2d}{}{\[\displaystyle \]}
\begin{Maple Normal}{
\begin{Maple Normal}{
}\end{Maple Normal}
}\end{Maple Normal}
\begin{Maple Normal}{
\begin{Maple Normal}{
\mapleinline{inert}{2d}{}{\[\displaystyle \]}
}\end{Maple Normal}
}\end{Maple Normal}

%% file: contas2.tex
\pagestyle{empty}
\DefineParaStyle{Maple Heading 1}
\DefineParaStyle{Maple Text Output}
\DefineParaStyle{Maple Dash Item}
\DefineParaStyle{Maple Bullet Item}
\DefineParaStyle{Maple Normal}
\DefineParaStyle{Maple Heading 4}
\DefineParaStyle{Maple Heading 3}
\DefineParaStyle{Maple Heading 2}
\DefineParaStyle{Maple Warning}
\DefineParaStyle{Maple Title}
\DefineParaStyle{Maple Error}
\DefineCharStyle{Maple Hyperlink}
\DefineCharStyle{Maple 2D Math}
\DefineCharStyle{Maple Maple Input}
\DefineCharStyle{Maple 2D Output}
\DefineCharStyle{Maple 2D Input}
\begin{maplegroup}
\begin{mapleinput}
\mapleinline{active}{2d}{libname := `C:/Finsler`, libname; 1}{\[\]}
\end{mapleinput}
\mapleresult
\begin{maplelatex}
\mapleinline{inert}{2d}{libname := "C:/Finsler", "C:\Program Files\Maple 17\lib", "."}{\[\displaystyle {\it libname}\, := \,``C:/Finsler'',\,``C:Program FilesMaple 17lib'',\,``.''\]}
\end{maplelatex}
\end{maplegroup}
\begin{maplegroup}
\begin{mapleinput}
\mapleinline{active}{2d}{with(Finsler); 1}{\[\]}
\end{mapleinput}
\mapleresult
\begin{maplelatex}
\mapleinline{inert}{2d}{[Dcoordinates, Hdiff, K, connection, init, metricfunction, tddiff]}{\[\displaystyle [{\it Dcoordinates},{\it Hdiff},K,{\it connection}\\
\mbox{},{\it init},{\it metricfunction},{\it tddiff}]\]}
\end{maplelatex}
\end{maplegroup}
\begin{maplegroup}
\begin{mapleinput}
\mapleinline{active}{2d}{dimension := 2; 1}{\[\]}
\end{mapleinput}
\mapleresult
\begin{maplelatex}
\mapleinline{inert}{2d}{dimension := 2}{\[\displaystyle {\it dimension}\, := \,2\]}
\end{maplelatex}
\end{maplegroup}
\begin{maplegroup}
\begin{mapleinput}
\mapleinline{active}{2d}{coordinates(x1, x2); 1}{\[\]}
\end{mapleinput}
\mapleresult
\begin{maplelatex}
\mapleinline{inert}{2d}{The coordinates are:}{\[\displaystyle \mbox {{\tt The coordinates are:}}\]}
\end{maplelatex}
\mapleresult
\begin{maplelatex}
\mapleinline{inert}{2d}{X* ^1 = x1}{\[\displaystyle X{\mbox {{\tt }}}^{\mbox {{\tt 1}}}={\it x1}\]}
\end{maplelatex}
\mapleresult
\begin{maplelatex}
\mapleinline{inert}{2d}{X*^2 = x2}{\[\displaystyle X{\mbox {{\tt }}}^{2}={\it x2}\]}
\end{maplelatex}
\end{maplegroup}
\begin{maplegroup}
\begin{mapleinput}
\mapleinline{active}{2d}{Dcoordinates(y1, y2); 1}{\[\]}
\end{mapleinput}
\mapleresult
Y assigned to DCoordinateName
\mapleresult
\begin{maplelatex}
\end{maplelatex}
\mapleresult
\begin{maplelatex}
\mapleinline{inert}{2d}{Y*^1 = y1}{\[\displaystyle Y{\mbox {{\tt }}}^{\mbox {{\tt 1}}}={\it y1}\]}
\end{maplelatex}
\mapleresult
\begin{maplelatex}
\mapleinline{inert}{2d}{Y*^2 = y2}{\[\displaystyle Y{\mbox {{\tt }}}^{2}={\it y2}\]}
\end{maplelatex}
\end{maplegroup}
\begin{maplegroup}
\begin{mapleinput}
\mapleinline{active}{2d}{F := exp(\nu_3*x1*x2-a*x1+b*x2)*y2^(1+1/lambda)/y1^(1/lambda); 1}{\[\]}
\end{mapleinput}
\mapleresult
\begin{maplelatex}
\mapleinline{inert}{2d}{Finsler Metric:}{\[\displaystyle \mbox {{\tt Finsler Metric:}}\]}
\end{maplelatex}
\begin{maplelatex}
\mapleinline{inert}{2d}{F := exp(\nu_3*x1*x2-a*x1+b*x2)*y2^(1+1/lambda)/y1^(1/lambda)}{\[\displaystyle F\, := \,{\frac {{{\rm e}^{\nu_3{\it x1}\,{\it x2}-a{\it x1}+b{\it x2}}}{{\it y2}}^{1+\frac{1}{\lambda}}}{{{\it y1}\\
\mbox{}}^{\frac{1}{\lambda}}}}\]}
\end{maplelatex}
\end{maplegroup}
\begin{maplegroup}
\begin{mapleinput}
\mapleinline{active}{2d}{}{\[\]}
\end{mapleinput}
\end{maplegroup}
\mapleinline{inert}{2d}{metricfunction(F^2); 1}{\[\displaystyle \]}
\begin{maplegroup}
\mapleresult
\begin{maplelatex}
\mapleinline{inert}{2d}{The components of the metric are:}{\[\displaystyle \mbox {{\tt The components of the metric are:}}\]}
\end{maplelatex}
\mapleresult
\begin{maplelatex}
\mapleinline{inert}{2d}{Typesetting:-mrow(Typesetting:-mrow(Typesetting:-mi("g", italic = "true", mathvariant = "italic"), Typesetting:-mo("&InvisibleTimes;", mathvariant = "normal", fence = "false", separator = "false", stretchy = "false", symmetric = "false", largeop = "false", movablelimits = "false", accent = "false", lspace = "0.0em", rspace = "0.0em"), Typesetting:-msub(Typesetting:-mi(""), Typesetting:-mrow(Typesetting:-mi("x1", italic = "true", mathvariant = "italic")), subscriptshift = "0"), Typesetting:-mo("&InvisibleTimes;", mathvariant = "normal", fence = "false", separator = "false", stretchy = "false", symmetric = "false", largeop = "false", movablelimits = "false", accent = "false", lspace = "0.0em", rspace = "0.0em"), Typesetting:-msub(Typesetting:-mi(""), Typesetting:-mrow(Typesetting:-mi("x1", italic = "true", mathvariant = "italic")), subscriptshift = "0")), Typesetting:-mo("=", mathvariant = "normal", fence = "false", separator = "false", stretchy = "false", symmetric = "false", largeop = "false", movablelimits = "false", accent = "false", lspace = "0.2777778em", rspace = "0.2777778em"), Typesetting:-mfrac(Typesetting:-mrow(Typesetting:-msup(Typesetting:-mfenced(Typesetting:-msup(Typesetting:-mo("&ExponentialE;", mathvariant = "normal", fence = "false", separator = "false", stretchy = "false", symmetric = "false", largeop = "false", movablelimits = "false", accent = "false", lspace = "0.0em", rspace = "0.1111111em"), Typesetting:-mrow(Typesetting:-mrow(Typesetting:-mi("c", italic = "true", mathvariant = "italic"), Typesetting:-mo("&InvisibleTimes;", mathvariant = "normal", fence = "false", separator = "false", stretchy = "false", symmetric = "false", largeop = "false", movablelimits = "false", accent = "false", lspace = "0.0em", rspace = "0.0em"), Typesetting:-mi("x1", italic = "true", mathvariant = "italic"), Typesetting:-mo("&InvisibleTimes;", mathvariant = "normal", fence = "false", separator = "false", stretchy = "false", symmetric = "false", largeop = "false", movablelimits = "false", accent = "false", lspace = "0.0em", rspace = "0.0em"), Typesetting:-mi("x2", italic = "true", mathvariant = "italic")), Typesetting:-mo("&minus;", mathvariant = "normal", fence = "false", separator = "false", stretchy = "false", symmetric = "false", largeop = "false", movablelimits = "false", accent = "false", lspace = "0.2222222em", rspace = "0.2222222em"), Typesetting:-mrow(Typesetting:-mi("a", italic = "true", mathvariant = "italic"), Typesetting:-mo("&InvisibleTimes;", mathvariant = "normal", fence = "false", separator = "false", stretchy = "false", symmetric = "false", largeop = "false", movablelimits = "false", accent = "false", lspace = "0.0em", rspace = "0.0em"), Typesetting:-mi("x1", italic = "true", mathvariant = "italic")), Typesetting:-mo("+", mathvariant = "normal", fence = "false", separator = "false", stretchy = "false", symmetric = "false", largeop = "false", movablelimits = "false", accent = "false", lspace = "0.2222222em", rspace = "0.2222222em"), Typesetting:-mrow(Typesetting:-mi("b", italic = "true", mathvariant = "italic"), Typesetting:-mo("&InvisibleTimes;", mathvariant = "normal", fence = "false", separator = "false", stretchy = "false", symmetric = "false", largeop = "false", movablelimits = "false", accent = "false", lspace = "0.0em", rspace = "0.0em"), Typesetting:-mi("x2", italic = "true", mathvariant = "italic"))), superscriptshift = "0"), mathvariant = "normal"), Typesetting:-mn("2", mathvariant = "normal"), superscriptshift = "0"), Typesetting:-mo("&InvisibleTimes;", mathvariant = "normal", fence = "false", separator = "false", stretchy = "false", symmetric = "false", largeop = "false", movablelimits = "false", accent = "false", lspace = "0.0em", rspace = "0.0em"), Typesetting:-msup(Typesetting:-mfenced(Typesetting:-msup(Typesetting:-mi("y2", italic = "true", mathvariant = "italic"), Typesetting:-mfrac(Typesetting:-mrow(Typesetting:-mrow(Typesetting:-mi("&lambda;", italic = "false", mathvariant = "normal"), Typesetting:-mo("+", mathvariant = "normal", fence = "false", separator = "false", stretchy = "false", symmetric = "false", largeop = "false", movablelimits = "false", accent = "false", lspace = "0.2222222em", rspace = "0.2222222em"), Typesetting:-mn("1", mathvariant = "normal"))), Typesetting:-mrow(Typesetting:-mi("&lambda;", italic = "false", mathvariant = "normal")), linethickness = "1", denomalign = "center", numalign = "center", bevelled = "false"), superscriptshift = "0"), mathvariant = "normal"), Typesetting:-mn("2", mathvariant = "normal"), superscriptshift = "0"), Typesetting:-mo("&InvisibleTimes;", mathvariant = "normal", fence = "false", separator = "false", stretchy = "false", symmetric = "false", largeop = "false", movablelimits = "false", accent = "false", lspace = "0.0em", rspace = "0.0em"), Typesetting:-mfenced(Typesetting:-mrow(Typesetting:-mi("&lambda;", italic = "false", mathvariant = "normal"), Typesetting:-mo("+", mathvariant = "normal", fence = "false", separator = "false", stretchy = "false", symmetric = "false", largeop = "false", movablelimits = "false", accent = "false", lspace = "0.2222222em", rspace = "0.2222222em"), Typesetting:-mn("2", mathvariant = "normal")), mathvariant = "normal")), Typesetting:-mrow(Typesetting:-msup(Typesetting:-mfenced(Typesetting:-msup(Typesetting:-mi("y1", italic = "true", mathvariant = "italic"), Typesetting:-mfrac(Typesetting:-mn("1", mathvariant = "normal"), Typesetting:-mrow(Typesetting:-mi("&lambda;", italic = "false", mathvariant = "normal")), linethickness = "1", denomalign = "center", numalign = "center", bevelled = "false"), superscriptshift = "0"), mathvariant = "normal"), Typesetting:-mn("2", mathvariant = "normal"), superscriptshift = "0"), Typesetting:-mo("&InvisibleTimes;", mathvariant = "normal", fence = "false", separator = "false", stretchy = "false", symmetric = "false", largeop = "false", movablelimits = "false", accent = "false", lspace = "0.0em", rspace = "0.0em"), Typesetting:-msup(Typesetting:-mi("&lambda;", italic = "false", mathvariant = "normal"), Typesetting:-mn("2", mathvariant = "normal"), superscriptshift = "0"), Typesetting:-mo("&InvisibleTimes;", mathvariant = "normal", fence = "false", separator = "false", stretchy = "false", symmetric = "false", largeop = "false", movablelimits = "false", accent = "false", lspace = "0.0em", rspace = "0.0em"), Typesetting:-msup(Typesetting:-mi("y1", italic = "true", mathvariant = "italic"), Typesetting:-mn("2", mathvariant = "normal"), superscriptshift = "0")), linethickness = "1", denomalign = "center", numalign = "center", bevelled = "false"))}{\[\displaystyle g ~_{{\it x1} }~_{{\it x1} }=\frac{\left(e^{\nu_3 ~{\it x1} ~{\it x2} -a ~{\it x1} +b ~{\it x2} }\right)^{2}~\left({\it y2} ^{\frac{\lambda +1}{\lambda }}\right)^{2}~\left(\lambda +2\right)}{\left({\it y1} ^{\frac{1}{\lambda }}\right)^{2}~\lambda ^{2}~{\it y1} ^{2}}\]}
\end{maplelatex}
\mapleresult
\begin{maplelatex}
\mapleinline{inert}{2d}{Typesetting:-mrow(Typesetting:-mrow(Typesetting:-mi("g", italic = "true", mathvariant = "italic"), Typesetting:-mo("&InvisibleTimes;", mathvariant = "normal", fence = "false", separator = "false", stretchy = "false", symmetric = "false", largeop = "false", movablelimits = "false", accent = "false", lspace = "0.0em", rspace = "0.0em"), Typesetting:-msub(Typesetting:-mi(""), Typesetting:-mrow(Typesetting:-mi("x1", italic = "true", mathvariant = "italic")), subscriptshift = "0"), Typesetting:-mo("&InvisibleTimes;", mathvariant = "normal", fence = "false", separator = "false", stretchy = "false", symmetric = "false", largeop = "false", movablelimits = "false", accent = "false", lspace = "0.0em", rspace = "0.0em"), Typesetting:-msub(Typesetting:-mi(""), Typesetting:-mrow(Typesetting:-mi("x2", italic = "true", mathvariant = "italic")), subscriptshift = "0")), Typesetting:-mo("=", mathvariant = "normal", fence = "false", separator = "false", stretchy = "false", symmetric = "false", largeop = "false", movablelimits = "false", accent = "false", lspace = "0.2777778em", rspace = "0.2777778em"), Typesetting:-mrow(Typesetting:-mo("&uminus0;", mathvariant = "normal", fence = "false", separator = "false", stretchy = "false", symmetric = "false", largeop = "false", movablelimits = "false", accent = "false", lspace = "0.2222222em", rspace = "0.2222222em"), Typesetting:-mfrac(Typesetting:-mrow(Typesetting:-mn("2", mathvariant = "normal"), Typesetting:-mo("&InvisibleTimes;", mathvariant = "normal", fence = "false", separator = "false", stretchy = "false", symmetric = "false", largeop = "false", movablelimits = "false", accent = "false", lspace = "0.0em", rspace = "0.0em"), Typesetting:-msup(Typesetting:-mfenced(Typesetting:-msup(Typesetting:-mo("&ExponentialE;", mathvariant = "normal", fence = "false", separator = "false", stretchy = "false", symmetric = "false", largeop = "false", movablelimits = "false", accent = "false", lspace = "0.0em", rspace = "0.1111111em"), Typesetting:-mrow(Typesetting:-mrow(Typesetting:-mi("c", italic = "true", mathvariant = "italic"), Typesetting:-mo("&InvisibleTimes;", mathvariant = "normal", fence = "false", separator = "false", stretchy = "false", symmetric = "false", largeop = "false", movablelimits = "false", accent = "false", lspace = "0.0em", rspace = "0.0em"), Typesetting:-mi("x1", italic = "true", mathvariant = "italic"), Typesetting:-mo("&InvisibleTimes;", mathvariant = "normal", fence = "false", separator = "false", stretchy = "false", symmetric = "false", largeop = "false", movablelimits = "false", accent = "false", lspace = "0.0em", rspace = "0.0em"), Typesetting:-mi("x2", italic = "true", mathvariant = "italic")), Typesetting:-mo("&minus;", mathvariant = "normal", fence = "false", separator = "false", stretchy = "false", symmetric = "false", largeop = "false", movablelimits = "false", accent = "false", lspace = "0.2222222em", rspace = "0.2222222em"), Typesetting:-mrow(Typesetting:-mi("a", italic = "true", mathvariant = "italic"), Typesetting:-mo("&InvisibleTimes;", mathvariant = "normal", fence = "false", separator = "false", stretchy = "false", symmetric = "false", largeop = "false", movablelimits = "false", accent = "false", lspace = "0.0em", rspace = "0.0em"), Typesetting:-mi("x1", italic = "true", mathvariant = "italic")), Typesetting:-mo("+", mathvariant = "normal", fence = "false", separator = "false", stretchy = "false", symmetric = "false", largeop = "false", movablelimits = "false", accent = "false", lspace = "0.2222222em", rspace = "0.2222222em"), Typesetting:-mrow(Typesetting:-mi("b", italic = "true", mathvariant = "italic"), Typesetting:-mo("&InvisibleTimes;", mathvariant = "normal", fence = "false", separator = "false", stretchy = "false", symmetric = "false", largeop = "false", movablelimits = "false", accent = "false", lspace = "0.0em", rspace = "0.0em"), Typesetting:-mi("x2", italic = "true", mathvariant = "italic"))), superscriptshift = "0"), mathvariant = "normal"), Typesetting:-mn("2", mathvariant = "normal"), superscriptshift = "0"), Typesetting:-mo("&InvisibleTimes;", mathvariant = "normal", fence = "false", separator = "false", stretchy = "false", symmetric = "false", largeop = "false", movablelimits = "false", accent = "false", lspace = "0.0em", rspace = "0.0em"), Typesetting:-msup(Typesetting:-mfenced(Typesetting:-msup(Typesetting:-mi("y2", italic = "true", mathvariant = "italic"), Typesetting:-mfrac(Typesetting:-mrow(Typesetting:-mrow(Typesetting:-mi("&lambda;", italic = "false", mathvariant = "normal"), Typesetting:-mo("+", mathvariant = "normal", fence = "false", separator = "false", stretchy = "false", symmetric = "false", largeop = "false", movablelimits = "false", accent = "false", lspace = "0.2222222em", rspace = "0.2222222em"), Typesetting:-mn("1", mathvariant = "normal"))), Typesetting:-mrow(Typesetting:-mi("&lambda;", italic = "false", mathvariant = "normal")), linethickness = "1", denomalign = "center", numalign = "center", bevelled = "false"), superscriptshift = "0"), mathvariant = "normal"), Typesetting:-mn("2", mathvariant = "normal"), superscriptshift = "0"), Typesetting:-mo("&InvisibleTimes;", mathvariant = "normal", fence = "false", separator = "false", stretchy = "false", symmetric = "false", largeop = "false", movablelimits = "false", accent = "false", lspace = "0.0em", rspace = "0.0em"), Typesetting:-mfenced(Typesetting:-mrow(Typesetting:-mi("&lambda;", italic = "false", mathvariant = "normal"), Typesetting:-mo("+", mathvariant = "normal", fence = "false", separator = "false", stretchy = "false", symmetric = "false", largeop = "false", movablelimits = "false", accent = "false", lspace = "0.2222222em", rspace = "0.2222222em"), Typesetting:-mn("1", mathvariant = "normal")), mathvariant = "normal")), Typesetting:-mrow(Typesetting:-msup(Typesetting:-mfenced(Typesetting:-msup(Typesetting:-mi("y1", italic = "true", mathvariant = "italic"), Typesetting:-mfrac(Typesetting:-mn("1", mathvariant = "normal"), Typesetting:-mrow(Typesetting:-mi("&lambda;", italic = "false", mathvariant = "normal")), linethickness = "1", denomalign = "center", numalign = "center", bevelled = "false"), superscriptshift = "0"), mathvariant = "normal"), Typesetting:-mn("2", mathvariant = "normal"), superscriptshift = "0"), Typesetting:-mo("&InvisibleTimes;", mathvariant = "normal", fence = "false", separator = "false", stretchy = "false", symmetric = "false", largeop = "false", movablelimits = "false", accent = "false", lspace = "0.0em", rspace = "0.0em"), Typesetting:-msup(Typesetting:-mi("&lambda;", italic = "false", mathvariant = "normal"), Typesetting:-mn("2", mathvariant = "normal"), superscriptshift = "0"), Typesetting:-mo("&InvisibleTimes;", mathvariant = "normal", fence = "false", separator = "false", stretchy = "false", symmetric = "false", largeop = "false", movablelimits = "false", accent = "false", lspace = "0.0em", rspace = "0.0em"), Typesetting:-mi("y1", italic = "true", mathvariant = "italic"), Typesetting:-mo("&InvisibleTimes;", mathvariant = "normal", fence = "false", separator = "false", stretchy = "false", symmetric = "false", largeop = "false", movablelimits = "false", accent = "false", lspace = "0.0em", rspace = "0.0em"), Typesetting:-mi("y2", italic = "true", mathvariant = "italic")), linethickness = "1", denomalign = "center", numalign = "center", bevelled = "false")))}{\[\displaystyle g ~_{{\it x1} }~_{{\it x2} }=-\frac{2~\left(e^{\nu_3 ~{\it x1} ~{\it x2} -a ~{\it x1} +b ~{\it x2} }\right)^{2}~\left({\it y2} ^{\frac{\lambda +1}{\lambda }}\right)^{2}~\left(\lambda +1\right)}{\left({\it y1} ^{\frac{1}{\lambda }}\right)^{2}~\lambda ^{2}~{\it y1} ~{\it y2} }\]}
\end{maplelatex}
\mapleresult
\mapleresult
\begin{maplelatex}
\mapleinline{inert}{2d}{Typesetting:-mrow(Typesetting:-mrow(Typesetting:-mi("g", italic = "true", mathvariant = "italic"), Typesetting:-mo("&InvisibleTimes;", mathvariant = "normal", fence = "false", separator = "false", stretchy = "false", symmetric = "false", largeop = "false", movablelimits = "false", accent = "false", lspace = "0.0em", rspace = "0.0em"), Typesetting:-msub(Typesetting:-mi(""), Typesetting:-mrow(Typesetting:-mi("x2", italic = "true", mathvariant = "italic")), subscriptshift = "0"), Typesetting:-mo("&InvisibleTimes;", mathvariant = "normal", fence = "false", separator = "false", stretchy = "false", symmetric = "false", largeop = "false", movablelimits = "false", accent = "false", lspace = "0.0em", rspace = "0.0em"), Typesetting:-msub(Typesetting:-mi(""), Typesetting:-mrow(Typesetting:-mi("x2", italic = "true", mathvariant = "italic")), subscriptshift = "0")), Typesetting:-mo("=", mathvariant = "normal", fence = "false", separator = "false", stretchy = "false", symmetric = "false", largeop = "false", movablelimits = "false", accent = "false", lspace = "0.2777778em", rspace = "0.2777778em"), Typesetting:-mfrac(Typesetting:-mrow(Typesetting:-msup(Typesetting:-mfenced(Typesetting:-msup(Typesetting:-mo("&ExponentialE;", mathvariant = "normal", fence = "false", separator = "false", stretchy = "false", symmetric = "false", largeop = "false", movablelimits = "false", accent = "false", lspace = "0.0em", rspace = "0.1111111em"), Typesetting:-mrow(Typesetting:-mrow(Typesetting:-mi("c", italic = "true", mathvariant = "italic"), Typesetting:-mo("&InvisibleTimes;", mathvariant = "normal", fence = "false", separator = "false", stretchy = "false", symmetric = "false", largeop = "false", movablelimits = "false", accent = "false", lspace = "0.0em", rspace = "0.0em"), Typesetting:-mi("x1", italic = "true", mathvariant = "italic"), Typesetting:-mo("&InvisibleTimes;", mathvariant = "normal", fence = "false", separator = "false", stretchy = "false", symmetric = "false", largeop = "false", movablelimits = "false", accent = "false", lspace = "0.0em", rspace = "0.0em"), Typesetting:-mi("x2", italic = "true", mathvariant = "italic")), Typesetting:-mo("&minus;", mathvariant = "normal", fence = "false", separator = "false", stretchy = "false", symmetric = "false", largeop = "false", movablelimits = "false", accent = "false", lspace = "0.2222222em", rspace = "0.2222222em"), Typesetting:-mrow(Typesetting:-mi("a", italic = "true", mathvariant = "italic"), Typesetting:-mo("&InvisibleTimes;", mathvariant = "normal", fence = "false", separator = "false", stretchy = "false", symmetric = "false", largeop = "false", movablelimits = "false", accent = "false", lspace = "0.0em", rspace = "0.0em"), Typesetting:-mi("x1", italic = "true", mathvariant = "italic")), Typesetting:-mo("+", mathvariant = "normal", fence = "false", separator = "false", stretchy = "false", symmetric = "false", largeop = "false", movablelimits = "false", accent = "false", lspace = "0.2222222em", rspace = "0.2222222em"), Typesetting:-mrow(Typesetting:-mi("b", italic = "true", mathvariant = "italic"), Typesetting:-mo("&InvisibleTimes;", mathvariant = "normal", fence = "false", separator = "false", stretchy = "false", symmetric = "false", largeop = "false", movablelimits = "false", accent = "false", lspace = "0.0em", rspace = "0.0em"), Typesetting:-mi("x2", italic = "true", mathvariant = "italic"))), superscriptshift = "0"), mathvariant = "normal"), Typesetting:-mn("2", mathvariant = "normal"), superscriptshift = "0"), Typesetting:-mo("&InvisibleTimes;", mathvariant = "normal", fence = "false", separator = "false", stretchy = "false", symmetric = "false", largeop = "false", movablelimits = "false", accent = "false", lspace = "0.0em", rspace = "0.0em"), Typesetting:-msup(Typesetting:-mfenced(Typesetting:-msup(Typesetting:-mi("y2", italic = "true", mathvariant = "italic"), Typesetting:-mfrac(Typesetting:-mrow(Typesetting:-mrow(Typesetting:-mi("&lambda;", italic = "false", mathvariant = "normal"), Typesetting:-mo("+", mathvariant = "normal", fence = "false", separator = "false", stretchy = "false", symmetric = "false", largeop = "false", movablelimits = "false", accent = "false", lspace = "0.2222222em", rspace = "0.2222222em"), Typesetting:-mn("1", mathvariant = "normal"))), Typesetting:-mrow(Typesetting:-mi("&lambda;", italic = "false", mathvariant = "normal")), linethickness = "1", denomalign = "center", numalign = "center", bevelled = "false"), superscriptshift = "0"), mathvariant = "normal"), Typesetting:-mn("2", mathvariant = "normal"), superscriptshift = "0"), Typesetting:-mo("&InvisibleTimes;", mathvariant = "normal", fence = "false", separator = "false", stretchy = "false", symmetric = "false", largeop = "false", movablelimits = "false", accent = "false", lspace = "0.0em", rspace = "0.0em"), Typesetting:-mfenced(Typesetting:-mrow(Typesetting:-mi("&lambda;", italic = "false", mathvariant = "normal"), Typesetting:-mo("+", mathvariant = "normal", fence = "false", separator = "false", stretchy = "false", symmetric = "false", largeop = "false", movablelimits = "false", accent = "false", lspace = "0.2222222em", rspace = "0.2222222em"), Typesetting:-mn("1", mathvariant = "normal")), mathvariant = "normal"), Typesetting:-mo("&InvisibleTimes;", mathvariant = "normal", fence = "false", separator = "false", stretchy = "false", symmetric = "false", largeop = "false", movablelimits = "false", accent = "false", lspace = "0.0em", rspace = "0.0em"), Typesetting:-mfenced(Typesetting:-mrow(Typesetting:-mi("&lambda;", italic = "false", mathvariant = "normal"), Typesetting:-mo("+", mathvariant = "normal", fence = "false", separator = "false", stretchy = "false", symmetric = "false", largeop = "false", movablelimits = "false", accent = "false", lspace = "0.2222222em", rspace = "0.2222222em"), Typesetting:-mn("2", mathvariant = "normal")), mathvariant = "normal")), Typesetting:-mrow(Typesetting:-msup(Typesetting:-mfenced(Typesetting:-msup(Typesetting:-mi("y1", italic = "true", mathvariant = "italic"), Typesetting:-mfrac(Typesetting:-mn("1", mathvariant = "normal"), Typesetting:-mrow(Typesetting:-mi("&lambda;", italic = "false", mathvariant = "normal")), linethickness = "1", denomalign = "center", numalign = "center", bevelled = "false"), superscriptshift = "0"), mathvariant = "normal"), Typesetting:-mn("2", mathvariant = "normal"), superscriptshift = "0"), Typesetting:-mo("&InvisibleTimes;", mathvariant = "normal", fence = "false", separator = "false", stretchy = "false", symmetric = "false", largeop = "false", movablelimits = "false", accent = "false", lspace = "0.0em", rspace = "0.0em"), Typesetting:-msup(Typesetting:-mi("&lambda;", italic = "false", mathvariant = "normal"), Typesetting:-mn("2", mathvariant = "normal"), superscriptshift = "0"), Typesetting:-mo("&InvisibleTimes;", mathvariant = "normal", fence = "false", separator = "false", stretchy = "false", symmetric = "false", largeop = "false", movablelimits = "false", accent = "false", lspace = "0.0em", rspace = "0.0em"), Typesetting:-msup(Typesetting:-mi("y2", italic = "true", mathvariant = "italic"), Typesetting:-mn("2", mathvariant = "normal"), superscriptshift = "0")), linethickness = "1", denomalign = "center", numalign = "center", bevelled = "false"))}{\[\displaystyle g ~_{{\it x2} }~_{{\it x2} }=\frac{\left(e^{\nu_3 ~{\it x1} ~{\it x2} -a ~{\it x1} +b ~{\it x2} }\right)^{2}~\left({\it y2} ^{\frac{\lambda +1}{\lambda }}\right)^{2}~\left(\lambda +1\right)~\left(\lambda +2\right)}{\left({\it y1} ^{\frac{1}{\lambda }}\right)^{2}~\lambda ^{2}~{\it y2} ^{2}}\]}
\end{maplelatex}
\end{maplegroup}
\mapleinline{inert}{2d}{show(G[i]); 1}{\[\displaystyle \]}
\begin{maplegroup}
\begin{mapleinput}
\mapleinline{active}{2d}{}{\[\]}
\end{mapleinput}
\end{maplegroup}
\begin{maplegroup}
\mapleresult

\begin{maplelatex}
\mapleinline{inert}{2d}{G*^x1 = (1/2)*y1^2*lambda*(-c*x2+a)}{\[\displaystyle G{\mbox {{\tt }}}^{{\it x1}}=\frac{1}{2}\,{{\it y1}}^{2}\lambda\, \left( -\nu_3{\it x2}+a \right) \]}
\end{maplelatex}
\mapleresult
\begin{maplelatex}
\mapleinline{inert}{2d}{G*^x2 = (1/2)*y2^2*lambda*(c*x1+b)/(lambda+1)}{\[\displaystyle G{\mbox {{\tt }}}^{{\it x2}}=\frac{1}{2}\,{\frac {{{\it y2}}^{2}\lambda\, \left( \nu_3{\it x1}+b \right) \\
\mbox{}}{\lambda+1}}\]}
\end{maplelatex}
\end{maplegroup}
\mapleinline{inert}{2d}{K := K(a, b); 1}{\[\displaystyle \]}
\begin{maplegroup}
\mapleresult
\begin{maplelatex}
\mapleinline{inert}{2d}{Curvature:}{\[\displaystyle \mbox {{\tt Curvature:}}\]}
\end{maplelatex}
\begin{maplelatex}
\mapleinline{inert}{2d}{y2*y1*(y1^(1/lambda))^2*lambda^2*\nu_3/((exp(c*x1*x2-a*x1+b*x2))^2*(y2^((lambda+1)/lambda))^2*(lambda+1))}{\[\mathcal{K}=\nu_3\displaystyle\frac{y2 y1\left(y1^{\frac{1}{\lambda}}\right)^2 \lambda^2 }{\left(\mathrm{e}^{\nu_3 x1 x 2-a x1+b x 2}\right)^2\left(y 2^{\frac{\lambda+1}{\lambda}}\right)^2(\lambda+1)}\]}
\end{maplelatex}
\end{maplegroup}
\begin{Maple Normal}{
\begin{Maple Normal}{
\mapleinline{inert}{2d}{}{\[\displaystyle \]}
}\end{Maple Normal}
}\end{Maple Normal}